\documentclass[a4paper,oneside]{article}
\usepackage{amsthm,amsmath,amssymb,color}
\newtheorem{teor}{Theorem}[section]
\newtheorem{defin}[teor]{Definition}
\newtheorem{lemm}[teor]{Lemma}
\newtheorem{osse}[teor]{Remark}
\newtheorem{prop}[teor]{Proposition}
\newtheorem{defi}[teor]{Definition}
\newtheorem{coro}[teor]{Corollary}
\newtheorem{prob}[teor]{Problem}

\newcommand{\bele}{\begin{lemm}\begin{sl}}
\newcommand{\enle}{\end{sl}\end{lemm}}
\newcommand{\bedef}{\begin{defi}\begin{sl}}
\newcommand{\eddef}{\end{sl}\end{defi}}
\newcommand{\bete}{\begin{teor}\begin{sl}}
\newcommand{\ente}{\end{sl}\end{teor}}
\newcommand{\beos}{\begin{osse}\begin{rm}}
\newcommand{\eddos}{\end{rm}\end{osse}}
\newcommand{\bepr}{\begin{prop}\begin{sl}}
\newcommand{\empr}{\end{sl}\end{prop}}
\newcommand{\bepro}{\begin{prob}\begin{rm}}
\newcommand{\empro}{\end{rm}\end{prob}}
\newcommand{\bede}{\begin{defin}\begin{sl}}
\newcommand{\edde}{\end{sl}\end{defin}}
\newcommand{\beco}{\begin{coro}\begin{sl}}
\newcommand{\enco}{\end{sl}\end{coro}}

\newcommand{\quand}{\quad\text{and}\quad}
\newcommand{\qquand}{\qquad\text{and}\qquad}
\newcommand{\quext}{\quad\text}
\newcommand{\qquext}{\qquad\text}

\newcommand{\Lbfi}[1]{L^{#1}_{b,\phi}}

\newcommand{\ltwob}{L^2_{b}(\Omega)}
\newcommand{\ltwoloc}{L^2_{\loc}(\Omega)}
\newcommand{\ltwox}{L^2_{\barx}(\Omega)}
\newcommand{\wtwob}{W^{1,2}_{b}(\Omega)}
\newcommand{\wtwox}{W^{1,2}_{\barx}(\Omega)}
\newcommand{\wmtwox}{W^{-1,2}_{\barx}(\Omega)}
\newcommand{\lpx}{L^p_{\barx}(\Omega)}
\newcommand{\lppx}{L^{p'}_{\barx}(\Omega)}
\newcommand{\wnegx}{W^{-1,2}_{\barx}(\Omega)}

\newcommand{\bnul}{\mathcal{B}}
\newcommand{\norm}[3]{\|{#1}\|_{#2}^{#3}}
\newcommand{\dual}[2]{\langle #1, #2 \rangle}

\newcommand{\ltwobfi}{L^2_{b,\phi}(\Omega)}

\newcommand{\WQbfi}{W_{b,\phi}(Q)}
\newcommand{\XQbfi}{\Lbfi{2}(0,\ell;L^2(\Omega))}
\newcommand{\XQbpsi}{L^2_{b,\Psixr}(0,\ell;L^2(\Omega))}
\newcommand{\VQbfi}{L^2_{b,\phi}(0,\ell;W^{1,2}(\Omega))}
\newcommand{\ZQbfi}{L^2_{b,\phi}(0,\ell;W^{-1,2}(\Omega))}

\newcommand{\dwxt}{e^{-|x-\barx|}\,\dix \dit}
\newcommand{\dwx}{e^{-|x-\barx|}\,\dix}
\newcommand{\iq}{\int_{\Omega\times(0,\ell)}}

\newcommand{\xxQ}{\mathcal{X}}
\newcommand{\atr}{\mathcal{A}}
\newcommand{\atrQ}{\atr_{\ell}}

\newcommand{\hent}[3]{H_{#1}(#2,#3)}
\newcommand{\oo}{\mathcal{O}}
\newcommand{\Indo}{\mathbb{I}(\oo)}
\newcommand{\Reps}{R(\varepsilon)}
\newcommand{\epsn}{\eps_0}
\newcommand{\Psixr}{\psi_{x_0,R}}
\newcommand{\Omxr}{\Omega_{x_0,R}}

\newcommand{\RR}{\mathbb{R}}

\newcommand{\NN}{\mathbb{N}}

\newcommand{\ZZ}{\mathbb{Z}}

\newcommand{\beeq}[1]{\begin{equation}\label{#1}}
\newcommand{\eddeq}{\end{equation}}

\newcommand{\beeqa}[1]{\begin{eqnarray}\label{#1}}
\newcommand{\eddeqa}{\end{eqnarray}}

\newcommand{\beal}[1]{\begin{align}\label{#1}}
\newcommand{\eddal}{\end{align}}

\newcommand{\bespl}[1]{\begin{split}\label{#1}}
\newcommand{\edspl}{\end{split}}

\newcommand{\bega}[1]{\begin{gather}\label{#1}}
\newcommand{\edga}{\end{gather}}

\newcommand{\beeqax}{\begin{eqnarray*}}
\newcommand{\eddeqax}{\end{eqnarray*}}

\def\qed{\ifmmode 
  \else \leavevmode\unskip\penalty9999 \hbox{}\nobreak\hfill
  \fi
  \quad\hbox{\hskip.5em\vrule width.4em height.6em depth.05em\hskip.1em}}
\def\endproofsym{\qed}
\renewenvironment{proof}[1][Proof]{\trivlist\item[\hskip\labelsep{\hskip0pt
    {\normalfont\scshape#1.}\hskip .321429\parindent}]\ignorespaces}
{\endproofsym\endtrivlist}
\def\endnobox{\def\endproofsym{}\end{proof}\def\endproofsym{\qed}}

\newcommand{\no}{\nonumber}

\newcommand{\beeqao}{\begin{eqnarray}\no}
\newcommand{\bealo}{\begin{align}\no}
\newcommand{\besplo}{\begin{split}\no}
\newcommand{\begao}{\begin{gather}\no}

\newcommand{\intck}{\int_{C_k}}

\newcommand{\nor}[2]{\|#1\|_{#2}}

\newcommand{\eps}{\varepsilon}

\newcommand{\duav}[1]{\langle{#1}\rangle}

\newcommand{\ena}[1]{\e^{#1}}

\newcommand{\dt}{\partial_t}

\newcommand{\perogni}{\forall\,}

\newcommand{\io}{\int_\Omega}

\newcommand{\iTT}{\int_0^T}

\newcommand{\iTo}{\iTT\!\io}

\newcommand{\aalpha}{\boldsymbol{\alpha}}

\newcommand{\lhs}{left-hand side}
\newcommand{\rhs}{right-hand side}

\DeclareMathOperator{\dive}{div}
\DeclareMathOperator{\vol}{vol}
\DeclareMathOperator{\supp}{supp}

\DeclareMathOperator{\deriv}{d}

\DeclareMathOperator{\loc}{loc}

\DeclareMathOperator{\e}{e}

\let\TeXchi\chi
\def\chi{{\setbox0 \hbox{\mathsurround0pt
$\TeXchi$}\hbox{\raise\dp0 \copy0 }}}

\newcommand{\zzn}{_{0,n}}

\newcommand{\calX}{{\mathcal X}}
\newcommand{\calI}{{\mathcal I}}

\newcommand{\calD}{{\mathcal D}}

\newcommand{\barx}{\overline{x}}

\newcommand{\dit}{\deriv\!t}
\newcommand{\dis}{\deriv\!s}
\newcommand{\dix}{\deriv\!x}

\newcommand{\ddt}{\frac{\deriv\!{}}{\dit}}

\newcommand{\h}{{h}}

\definecolor{ddmagenta}{rgb}{0.7,0,1.0}
\definecolor{ddcyan}{rgb}{0,0.1,1.0}
\definecolor{dred}{rgb}{.8,0,0}

\numberwithin{equation}{section}

\begin{document}

\title{Attractors for nonlinear reaction-diffusion systems 
  in unbounded domains via the method of short trajectories}

\author{Maurizio Grasselli\\
{\sl Dipartimento di Matematica}\\
{\sl Politecnico di Milano}\\
{\sl Via Bonardi 9}\\
{\sl I-20133 Milano, Italy}\\
{\rm E-mail:~~\tt maurizio.grasselli@polimi.it}
\and
Dalibor Pra\v z\'ak\\
{\sl Katedra matematick\'e anal\'yzy}\\
{\sl Matematicko-fyzik\'aln\'{\i} fakulta Univerzity Karlovy}\\
{\sl Sokolovsk\'a 83, CZ-186 75, Czech Republic}\\
{\rm E-mail:~~\tt prazak@karlin.mff.cuni.cz}
\and
Giulio Schimperna\\
{\sl Dipartimento di Matematica, Universit\`a di Pavia}\\
{\sl Via Ferrata, 1, I-27100 Pavia, Italy}\\
{\rm E-mail:~~\tt giusch04@unipv.it}
}
\maketitle

\begin{abstract}
\noindent
We consider a nonlinear reaction-diffusion equation on the whole space $\RR^d$.
We prove the well-posedness of the corresponding Cauchy problem in a general
functional setting, namely, when the initial datum is uniformly locally bounded in 
$L^2$. Then we adapt the short trajectory method to establish the existence of the
global attractor and, if $d\leq 3$, we find an upper bound of its
Kolmogorov's $\varepsilon$-entropy.
\end{abstract}

\medskip


\medskip

\noindent
{\bf AMS (MOS) subject clas\-si\-fi\-ca\-tion:}~~35B41, 35K57, 92D25.


\section{Introduction}
\label{secintro}
The asymptotic behavior of solutions to reaction-diffusion equations 
has been the object of a large number of investigations. In particular,
the existence of global and exponential attractors and their fractal
dimension have been carefully studied in the case of bounded domains
(see, e.g., \cite{MZ} and references therein). 
However, unbounded domains are also rather
interesting for applications. In this case, the dynamics can exhibit
a more complex behavior characterized, for instance, by travelling waves 
connecting constant equilibria or by a continuum of space periodic 
equilibria (and much more, as shown in \cite{Ze03}). The lack of standard compact injections
require new ideas and the choice of the topology becomes crucial in order
not to exclude interesting invariant solutions. 
In the pioneering papers on the existence of global attractors \cite{A,BaVi90}
weighted norms were introduced and used. Since then, many contributions have
followed (see, e.g., \cite{ACDR1,ACDR2,EM,EZ,FLST,Me,Pr,PR,Wa} and references therein) 
making use of weighted or not phase spaces and under various assumptions 
on the reaction (and, possibly, convective) terms. 
However, as noticed in \cite{BaVi90}, the global attractor
can be noncompact (but just locally compact) and infinite dimensional. 
Actually, this is the more realistic case where the richness of the dynamics is preserved.
However, it is still possible to give a quantitative estimate of the thickness
of the global attractor by means of the so-called Kolmogorov's 
$\varepsilon$-entropy. Estimates of this quantity were proven in \cite{Ze01} under
quite general assumptions (see also \cite{EZ2} for a generalization which accounts for convection
and \cite{Ze03} for a careful analysis of related spatially chaotic phenomena).
Here we use a somewhat simplified weighted space setting along the lines of
\cite{Ze01} to analyze a reaction-diffusion equation of the form
\begin{equation} \label{eqref}
  \dt u - \dive a(\nabla u) + f(u) + \h(\nabla u) = g, \qquad
  \text{in }\; \RR^d \times (0,+\infty),
\end{equation}
where $a,f,h$ are suitable nonlinear functions and $f$ has a polynomially controlled
growth. More precisely, we introduce and solve an appropriate weighted weak formulation of 
the Cauchy problem for \eqref{eqref} with $g$ and the initial datum uniformly locally bounded in $L^2$. 
Then, by adapting the short trajectory method (see \cite{MP}), we easily prove the 
existence of the global attractor in a $L^2_{loc}$-topology. Finally, using once more that
approach, we estimate its $\varepsilon$-entropy from above. The main novelty with respect to
the existing literature (and, in particular, to \cite{Ze01})
is that we can essentially work in the usual ``parabolic'' functional
setting. Thus we only need a handful of (relatively) simple estimates, regularity assumptions
on $a,f,h,g$ are very mild,  and the phase-space includes bounded functions (and more). 
In addition, (reasonable) nonlinear diffusion terms  
along with typical reaction terms of the form $f(u)=u^3 - \gamma u$, $\gamma>0$, can be handled easily.
Possible extensions to systems are also pointed out. Further extensions will include delay effects
(cf., e.g., \cite{GP} for bounded domains).

The plan of this paper goes as follows. The next Section~\ref{secspaces} is devoted to introduce the 
functional setup which is, of course, a bit more complicated than the one with bounded domains.
The notion of Kolmogorov's $\varepsilon$-entropy adapted to our framework is also introduced.
Well-posedness and regularity issues are analyzed in Section~\ref{secwell}. The existence
of the global attractor is proven in Section~\ref{secdynam} and an upper bound for
its $\varepsilon$-entropy is established in Section~\ref{entroest}.


\section{Functional spaces}
\label{secspaces}

There are three classes of function spaces to be used in this
paper. The standard Lebesgue and Sobolev spaces together with
their weighted variants are briefly recalled in Section~\ref{sec21}.
These spaces are mainly used for formulating the existence theorem.

Throughout the paper $c_1$, $c_2 \dots$ denote universal constants
whose meaning can change with the context, but which are 
independent on the data of the equation and also on the weight functions.
We also occasionally simplify the notation by writing $A\approx B$
meaning that $c_1 A \le B \le c_2 A$.

\par
In Section~\ref{sec22}, we introduce the so-called \emph{uniformly
bounded} spaces and provide some equivalent descriptions of their
norms. These spaces are aimed at describing the
dynamics associated with the equation, and formulating the main
results.
As we mentioned in the Introduction, here we mostly follow \cite{Ze01,Ze03}, though
in a slightly simplified setting (thanks to the fact that $\Omega=\RR^d$
has no boundary). 
We remark that we perform the analysis by taking $\Omega=\RR^d$ for
simplicity; however, obvious technical adjustments would allow us to
treat any suitably regular unbounded domain $\Omega$. 

\par
Finally, in Section~\ref{sec23}, we describe a class of spaces that
can be thought of as a parabolic version of the uniformly bounded
spaces. These spaces are the main technical novelty of the paper
and are also the crucial tool for the application of the ``method
of trajectories'' to the problem of the dimension of the attractor.

\subsection{Weighted Sobolev spaces}
\label{sec21}

For a domain $\oo\subset \RR^d$, we use $L^p(\oo)$, $W^{1,p}(\oo)$,
$W^{1,p}_0(\oo)$ and $W^{-1,p'}=[W^{1,p}_0(\oo)]'$ to denote the
standard Sobolev spaces. Observe that $W^{1,p}(\Omega) = W^{1,p}_0(\Omega)$
as $\Omega=\RR^d$ throughout the paper. By $L^p_{\loc}(\Omega)$, 
$W^{1,p}_{\loc}(\Omega)$ we denote their locally integrable variants.

A prominent role will be played by the weight functions $e^{-|\cdot - \barx|}$, 
$\barx \in \RR^d$. These give rise to weighted spaces $\lpx$, 
$W^{1,p}_{\barx}(\Omega)$, 
$W^{-1,p'}_{\barx}(\Omega) = [W^{1,p}_{\barx}(\Omega)]'$, 
given via the respective norms
\begin{align*}
\norm{u}{\lpx}{p} &= \io |u(x)|^p \dwx,
\\
\norm{u}{W^{1,p}_{\barx}(\Omega)}{p}
&= \io \big( |\nabla u(x)|^p + |u(x)|^p \big)\dwx, 
\\
\norm{u}{W^{-1,p'}_{\barx}(\Omega)}{}
&= \sup_{v} \io u(x)v(x) \dwx,
\end{align*}
the last supremum being taken over $v\in W^{1,p}_{\barx}(\Omega)$
with unit norm. These spaces share the usual good properties of
Sobolev spaces (separability, reflexivity), note also that
\begin{equation} \label{hold-o}
\lpx \subset L^q_{\barx}(\Omega),
\qquad p\ge q,
\end{equation}
since $\Omega$ has finite measure with respect to the weight $e^{-|\cdot - \barx|}$.
It is also easy to see that
the spaces $\lpx$, $L^p_{\overline{y}}(\Omega)$ in fact
coincide and the equivalence constants only depend on 
$|\barx - \overline{y}|$.

\par
Finally, we remark that $\ltwox$, $\wtwox$  are Hilbert spaces using
the obvious scalar product; however, it is worth noting that 
$\norm{\nabla \cdot}{\ltwox}{}$ is not an equivalent norm on $\wtwox$
(just consider a constant function). The notation 
$\langle\cdot,\cdot\rangle_{\barx}$ will stand for the duality
pairing between $W^{-1,2}_{\overline{x}}(\Omega)$ and $\wtwox$.

\subsection{Uniformly bounded spaces}
\label{sec22}

First, we introduce the space
$\ltwob$ of the {\sl uniformly locally
$L^2$-functions}\/ as 
\begin{equation}\label{defL2b}
  L^2_b(\Omega):=\big\{u\in L^2_{\loc}(\Omega):
   \sup_{x_0\in\Omega} \|u\|_{L^2(C(x_0))}<+\infty\big\}.  
\end{equation}  
Here and below in the paper, $C(x_0)$ denotes the closed unit
cube of $\RR^d$ centered at $x_0$, namely
$C(x_0)=\prod_{i=1,\dots,d}[x_{0,i}-1/2,x_{0,i}+1/2]$,
where $x_{0,i}$ are the components of $x_0$.
Clearly, $L^2_b(\Omega)$, endowed with the graph norm,
is a Banach space. An equivalent norm is given by
\begin{equation}\label{norL2b}
  \norm{u}{\ltwob}{} := \sup_{k} \norm{u}{L^2(C_k)}{}.
\end{equation}  
Here and in what follows, $C_k$ are an enumeration of
the unit cubes centered at $x_k\in (\ZZ/2)^d$.
An advantage of this norm is that the supremum
is taken over a {\sl countable}\/ family of cubes. 
Note also that, for later convenience,
we allow a partial superposition of the cubes.

We will also need the weighted analogue of 
$L^2_b(\Omega)$. Given $\mu\ge 0$, an 
{\sl admissible weight}\/ 
of rate of growth $\mu$ is a (measurable and bounded) function 
$\phi:\RR^N\to (0,+\infty)$ satisfying, 
for some $c\ge 1$,
\begin{equation} \label{weight-gr}
  c^{-1} e^{-\mu|x-y|}
   \le \phi(x) / \phi(y) \le c e^{\mu|x-y|},
   \qquad   \perogni x,y\in\RR^d,
\end{equation}
as well as the estimate 
\begin{equation} \label{weight-bla}
   |\nabla \phi(x)| \le |\phi(x)|.
\end{equation}
A typical example is given by the exponential 
$\phi(x) = e^{m|x-\barx|}$, where $\barx\in\RR^d$ 
and $m\in[-1,0]$, which of course has rate of growth 
$\mu=|m|$. In fact, we can observe 
that $|m|\le 1$ would be enough
in order to have \eqref{weight-gr}-\eqref{weight-bla}.
However, since we also need global boundedness
of $\phi$ in the sequel, we will only consider
negative exponential weights. 

It is easy to prove (see \cite[Prop.~1.3]{Ze01})
that if $\phi_1$ and $\phi_2$ are admissible 
weights of growth rates $\mu_1$ and $\mu_2$, then
$\max\{\phi_1,\phi_2\}$ and $\min\{\phi_1,\phi_2\}$ are
still admissible weights both having growth rate
$\max\{\mu_1,\mu_2\}$.

We now have the analogue of \eqref{defL2b},
i.e., the space of the functions which are uniformly
locally $L^2$ with respect to the weight $\phi$. This is 
defined as
\begin{equation}\label{defL2bphi}
  \ltwobfi:=\big\{u\in L^2_{\loc}(\Omega):
   \sup_{x_0\in\Omega} \phi^{1/2}(x_0) \|u\|_{L^2(C(x_0))}<+\infty\big\}.  
\end{equation}  
As before, we will take on $\ltwobfi$ the equivalent
norm
\begin{equation}\label{norL2bphi}
  \|u\|_{\ltwobfi}^{2}:=\sup_{k}
    \phi(x_k) \norm{u}{L^2(C_k)}{2}.
\end{equation}  
It is easy to check that $\ltwobfi$ is then a Banach
space. Note that the constant function $1$ is
an admissible weight with growth rate $0$ and
$L^2_{b,1}(\Omega)=L^2_{b}(\Omega)$. Moreover, if the weight $\phi$
is of the form $\phi(x)=e^{-\mu|x-\barx|}$ for some
$\barx\in\RR^d$ and for $\mu\in[0,1]$, then it 
is $L^2_{b}(\Omega) \subset \ltwobfi$ with continuous inclusion.

Given now an admissible weight $\phi$ with rate of growth
{\sl strictly}\/ smaller than 1, we also
define 
\begin{equation}\label{deftildeX}
  \tilde L^2_{b,\phi}(\Omega)
  :=\Big\{u\in L^2_{\loc}(\Omega):
    u \in \ltwox~\perogni \barx\in \Omega,~
    \sup_{\barx\in\Omega}\phi(\barx)\io 
    |u(x)|^2 e^{-|x-\barx|}\,\dix < +\infty \Big\}.
\end{equation}  
It is not difficult 
to prove that $\tilde L^2_{b,\phi}(\Omega)$,
endowed with the graph norm, is 
also a Banach space (note that it is not a Hilbert 
space). More precisely, we can prove
\bete\label{teoequiv1}
 The spaces $\ltwobfi$ and $\tilde L^2_{b,\phi}(\Omega)$
 coincide and, in particular, their norms are equivalent.
\ente
\begin{proof}
 Recall that $\{C_k\}_{k\in\NN}$ are an enumeration of the unit 
 cubes of $\Omega$ 
 centered in the points of $x_k\in (\ZZ/2)^d$. 
 It is then clear that, for fixed $\barx\in\Omega$, we have
 \begin{equation}\label{dalibor11}
   \phi(\barx) \io |u(x)|^2 e^{-|x-\barx|}\,\dix
    \le \sum_{k\in\NN} \intck \phi(\barx) |u(x)|^2 e^{-|x-\barx|}\,\dix.
 \end{equation}
 Let us also notice that, for $x\in C_k$, there hold
 \begin{equation}\label{dalibor11b}
   e^{-|x-\barx|} \le c_1 e^{-|x_k-\barx|} \qquand
   \phi(\barx) \le c_2 \phi(x_k) e^{\mu|x_k - \barx|},
 \end{equation}
 for suitable $c_1,c_2>0$ independent of $k,\barx$. 
 Hence,
 \begin{align}\no
   \sum_{k\in\NN} \intck \phi(\barx) |u(x)|^2 e^{-|x-\barx|}\,\dix
    & \le c_3 \sum_{k\in\NN} e^{-(1-\mu)|x_k -\barx|} \phi(x_k) \intck |u(x)|^2\,\dix\\
  \label{dalibor12}
    & \le c_3 \Big( \sup_{k\in\NN} \phi(x_k) \nor{u}{L^2(C_k)}^2 \Big)
     \sum_{k\in\NN} e^{-(1-\mu)|x_k -\barx|}.
 \end{align}
 Assuming $\mu < 1$, the sum is bounded independently of $\barx$.
 Passing to the supremum, this entails that
 \begin{equation}\label{stnor11}
    \|u\|_{\tilde L^2_{b,\phi}(\Omega)} \le c \|u\|_{\ltwobfi}.
 \end{equation}
 To prove the opposite inequality, note that, for $x\in C(\barx)$,
 \begin{equation}\label{dalibor13}
   1 \le c_4 e^{-|x-\barx|}.
 \end{equation}
 Hence,
 \begin{align}\no
   \phi(\barx) \|u\|_{L^2(C(\barx))}^2
    & = \phi(\barx) \int_{C(\barx)} |u(x)|^2 \,\dix\\
  \no
    & \le c_4 \phi(\barx) \int_{C(\barx)} |u(x)|^2 e^{-|x-\barx|}\,\dix\\
  \label{dalibor14}
    & \le c_4 \phi(\barx) \io |u(x)|^2 e^{-|x-\barx|}\,\dix.
 \end{align}
 The proof is complete.
\end{proof}
\noindent%
For $p\in[1,\infty)$, we can define analogously
as above the spaces $L^p_{b,\phi}(\Omega)$ and $\tilde L^p_{b,\phi}(\Omega)$,
where
\begin{align*}
\norm{u}{L^p_{b,\phi}(\Omega)}{p} &:= \sup_k \phi(x_k)\norm{u}{L^p(C_k)}{p},
\\
\norm{u}{\tilde L^p_{b,\phi}(\Omega)}{p}
&:= \sup_{\barx} \phi(\barx) \io |u(x)|^p \dwx.
\end{align*}
As above, one proves that $L^p_{b,\phi}(\Omega) = \tilde L^p_{b,\phi}(\Omega)$
provided that the growth rate of $\phi$ is smaller than $1$;
the equivalence relation can be succinctly written as
\begin{equation} \label{ekvx}
\norm{u}{L^p_{b,\phi}(\Omega)}{p}
\approx
\sup_{\barx} \phi(\barx) \norm{u}{\lpx}{p}.
\end{equation}
Note also that the equivalence constants only depend on $\mu$, $c$
in \eqref{weight-gr} and not on the particular expression of 
the weight function; this fact will
be used repeatedly in various a priori estimates.

\par
As a next step, we extend the above construction
to Sobolev spaces. First of all, 
given an admissible weight $\phi$,
$W^{1,2}_b(\Omega)$ and $W^{1,2}_{b,\phi}(\Omega)$ are 
defined as the spaces of $L^2_{\loc}(\Omega)$-functions
which belong to $\ltwob$ and, respectively, 
$\ltwobfi$ together with their first 
(partial, distributional)
derivatives. These are, of course, Banach spaces
with the natural norms modelled on \eqref{norL2b}
and \eqref{norL2bphi}. We can also define
$\tilde W^{1,2}_{b,\phi}(\Omega)$ as the space of 
$L^2_{\loc}(\Omega)$-functions such that
\begin{equation}\label{deftildeV}
  \sup_{\barx\in\Omega}\phi(\barx)\io 
   \big(|u(x)|^2+|\nabla u(x)|^2\big)
   e^{-|x-\barx|}\,\dix < +\infty.
\end{equation}  
In particular, we write $\tilde W^{1,2}_{b}(\Omega)$ in case
$\phi\equiv1$. The analogue of Theorem~\ref{teoequiv1},
whose proof is omitted for brevity since it does
not present further difficulties, then reads
\bete\label{teoequiv2}
 Given an admissible weight $\phi$ of growth
 rate $\mu\in[0,1)$, the spaces $\tilde W^{1,2}_{b,\phi}(\Omega)$ 
 and $W^{1,2}_{b,\phi}(\Omega)$
 coincide and their norms are equivalent.
\ente
\noindent%
Next, we come to {\sl negative order}\/ spaces.
Firstly, we define 
\begin{equation}\label{defZbphi}
  W^{-1,2}_{b,\phi}(\Omega):=\Big\{\zeta\in\calD'(\Omega)
  :\zeta|_{C(\barx)}\in W^{-1,2}(C(\barx))~\perogni\barx\in\Omega,
   ~\sup_{\barx\in\Omega}
      \phi(\barx)\|\zeta\|_{W^{-1,2}(C(\barx))}^2<+\infty\Big\}.
\end{equation}  
Of course, the above, endowed with the graph norm,
is a Banach space (and the supremum could be
restricted to $\barx\in(\ZZ/2)^d$, see the proof of 
the next theorem).
As before, we can also define
the counterpart $\tilde W^{-1,2}_{b,\phi}(\Omega)$.
Let us take first  
$u\in \tilde L^2_{b,\phi}(\Omega)$ and
set
\begin{equation}\label{norZbphi}
  \|u\|_{\tilde W^{-1,2}_{b,\phi}(\Omega)}:=\sup_{v}
    \sup_{\barx\in\Omega}
   \phi(\barx)\io u(x) v(x) e^{-|x-\barx|}\,\dix,
\end{equation}  
where the first supremum is taken with respect to
\begin{equation}\label{qualiv}
  \big\{v\in W^{1,2}_{b,\phi}(\Omega):
    \|v\|_{\tilde W^{1,2}_{b,\phi}(\Omega)}\le 1\big\}.
\end{equation}  
The space $\tilde W^{-1,2}_{b,\phi}(\Omega)$ is then defined as 
the completion of $\ltwobfi$ with respect to the
norm \eqref{norZbphi}.
\bete\label{teoequiv3}
 Given an admissible weight $\phi$ of growth
 rate $\mu\in[0,1)$, the spaces $\tilde W^{-1,2}_{b,\phi}(\Omega)$ 
 and $W^{-1,2}_{b,\phi}(\Omega)$ coincide and their norms are equivalent.
\ente
\begin{proof}
 Let $C_k$ and $x_k$ be as in the proof of Theorem~\ref{teoequiv1}.
 Take $u\in \ltwobfi$. Then, it is clear that
 \begin{equation}\label{dalibor21}
   \phi^{1/2}(x_k) \nor{u}{W^{-1,2}(C_k)}
    = \sup_v \phi^{1/2}(x_k) \intck u(x)v(x)\,\dix,
 \end{equation}
 where the supremum is referred to the 
 $v$'s in $W^{1,2}_0(C_k)$ with $\|v\|_{W^{1,2}_0(C_k)}\le 1$.
 Let us take any such $v$ and extend it 
 by zero outside $C_k$. Then,
 \begin{align}\no
   \phi^{1/2}(x_k) \intck u(x)v(x)\,\dix
    & = \phi^{1/2}(x_k) \io u(x) v(x) \,\dix\\
  \label{dalibor22}
    & = \phi(x_k) \io u(x)
     \underbrace{ \phi^{-1/2}(x_k) v(x) \ena{|x-x_k|} }_{\tilde{v}(x)}
       \ena{-|x-x_k|}\,\dix.
 \end{align}
 One easily verifies that $\tilde{v}$ (which is in fact only supported
 in $C_k$) belongs to $W^{1,2}_{b,\phi}(\Omega)$ and has 
 the norm smaller than
 some constant $c_1$. 
 Taking the suprema with respect to $v$ and $k$,
 we eventually get that
 \begin{equation}\label{dalibor23}
   \|u\|_{W^{-1,2}_{b,\phi}(\Omega)} 
   \le c \|u\|_{\tilde W^{-1,2}_{b,\phi}(\Omega)}.
 \end{equation}
 The proof of the opposite inequality is a little bit harder.
 Let $u\in \ltwobfi$, $v\in W^{1,2}_{b,\phi}(\Omega)$ and $\barx\in\Omega$. 
 Let also 
 $\{\psi_k\}_{k\in\NN}$ be a smooth partition of unity associated
 to the cubes $C_k$. Then,
 \begin{align}\no
   \phi(\barx)\io u(x) v(x) e^{-|x-\barx|}\,\dix
   & = \sum_{k\in\NN} \phi(\barx) \io u(x) (v\psi_k)(x)
        e^{-|x-\barx|}\,\dix\\
  \no
    & \le \sum_{k\in\NN} \phi^{1/2}(x_k) \intck u(x) (v\psi_k)(x)
       \frac{\phi(\barx)}{\phi^{1/2}(x_k)} e^{-|x-\barx|}\,\dix\\
  \no
    & \le \sum_{k\in\NN} \phi^{1/2}(x_k) \| u \|_{W^{-1,2}(C_k)} 
          \| V_k \|_{W^{1,2}_0(C_k)}\\
  \label{dalibor24}  
    & \le \|u\|_{W^{-1,2}_{b,\phi}(\Omega)}
         \sum_{k\in\NN} \| V_k \|_{W^{1,2}_0(C_k)},
 \end{align}
 where we have set 
 \begin{equation}\label{dalibor25}
   V_k(x):=v(x)\psi_k(x)\phi(\barx)\phi^{-1/2}(x_k)e^{-|x-\barx|}.
 \end{equation}
 Then, a direct computation (notice that
 the functions $\psi_k$ can be chosen
 uniformly bounded together with their first derivatives) 
 shows that 
 \begin{equation}\label{norVk}
   \| V_k \|_{W^{1,2}_0(C_k)} \le 
    c \phi(\barx)\phi^{-1/2}(x_k)e^{-|x_k-\barx|}
    \| v \|_{W^{1,2}(C_k)},
 \end{equation}  
 where $c$ is independent of $\barx,k$. Thus, coming back to
 \eqref{dalibor24} and using Theorem~\ref{teoequiv2} and 
 \eqref{dalibor11b} once more,
 we arrive at 
 \begin{align}\no
   \phi(\barx)\io u(x) v(x) e^{-|x-\barx|}\,\dix
    & \le c \|u\|_{W^{-1,2}_{b,\phi}(\Omega)}
         \sum_{k\in\NN} \Big( \phi(\barx)\phi^{-1}(x_k)e^{-|x_k-\barx|}
     \phi^{1/2}(x_k) \| v \|_{W^{1,2}(C_k)}\Big)\\
  \no
    & \le c \|u\|_{W^{-1,2}_{b,\phi}(\Omega)} 
    \| v \|_{W^{1,2}_{b,\phi}(\Omega)}
     \sum_{k\in\NN} e^{-(1-\mu)|x_k-\barx|}\\
  \label{dalibor26}     
    & \le c \|u\|_{W^{-1,2}_{b,\phi}(\Omega)} 
    \| v \|_{\tilde W^{1,2}_{b,\phi}(\Omega)}.
 \end{align}
 Thus, dividing by $\| v \|_{\tilde W^{1,2}_{b,\phi}(\Omega)}$ and taking
 the supremum first with respect to $\barx$ and
 then with respect to $v$, we obtain the opposite
 inequality of \eqref{dalibor23}.
  
 To conclude the proof, we observe that, a priori, the equivalence of
 the norms of $W^{-1,2}_{b,\phi}(\Omega)$ 
 and $\tilde W^{-1,2}_{b,\phi}(\Omega)$ has been proved
 just for the functions of $\ltwobfi$. However, it can be 
 easily extended to the whole spaces by means of a 
 standard density argument. 
\end{proof}

It is worth while observing that the above 
defined spaces share some usual properties
of Lebesgue spaces, as for example
\[
\norm{uv}{\Lbfi{r}(\Omega)}{}
\le
\norm{u}{\Lbfi{p}(\Omega)}{}
\norm{v}{\Lbfi{q}(\Omega)}{},
\qquad \frac1r = \frac1p + \frac1q.
\]
On the other hand, the Sobolev embedding $W^{1,2}_{b,\phi}(\Omega)
\subset \Lbfi{p}(\Omega)$, $p=2d/(d-2)$ does not hold (unless $\phi\equiv1$)
due to the incompatibility of the powers of $\phi(x_k)$.  
Also, the $\ltwobfi$ norm of $\nabla u$ is not an equivalent norm in 
$W^{1,2}_{b,\phi}(\Omega)$.

Finally, we  will need seminorms that correspond to restrictions
to some (bounded) subdomain $\oo \subset \Omega$.
For arbitrary $\mathcal{O}\subset \Omega$, we set
\begin{align*}
\Indo &:= \big\{ k\in \NN;\ C_k \cap \oo \neq \emptyset \big\},
\\
\norm{u}{L^2_{b,\phi}(\oo)}{2} &:= 
\sup_{k\in \Indo} 
\phi(x_k) \norm{u}{L^2(C_k)}{2}.
\end{align*}

\subsection{Parabolic uniformly bounded spaces}
\label{sec23}

As an auxiliary tool,
we will work with a sort of ``parabolic version'' of
uniformly local spaces -- a main technical novelty of the present paper.
This setup seems rather natural for the study of dynamics of parabolic-like
evolutionary problems in unbounded domains. 

\par
Given an admissible weight function $\phi$, we define
spaces $\XQbfi$, $\VQbfi$ 
and $\ZQbfi$, 
where for any function 
$u(x,t):\Omega \times (0,\ell) \to \RR$, we set
\begin{align*}
\norm{u}{\XQbfi}{} &= \sup_{k\in \NN} \phi^{1/2}(x_k)
            \norm{u}{L^2(0,\ell;L^2(C_k))}{},
            \\
\norm{u}{\VQbfi}{} 
&= \sup_{k\in \NN} \phi^{1/2}(x_k)
            \norm{u}{L^2(0,\ell;W^{1,2}(C_k))}{},
            \\
\norm{u}{\ZQbfi}{} 
&= \sup_{k\in \NN} \phi^{1/2}(x_k)
            \norm{u}{L^{2}(0,\ell;W^{-1,2}(C_k))}{}.
\end{align*}
We also introduce the space $\Lbfi{p}(0,\ell;L^p(\Omega))$ as
\[
\norm{u}{\Lbfi{p}(0,\ell;L^p(\Omega))}{} 
= \sup_{k\in \NN} \phi^{1/p}(x_k) 
            \norm{u}{L^p(0,\ell;L^p(C_k))}{}.
\]
As customary, we omit the symbol $\phi$ if $\phi\equiv1$.
We will also need localized seminorm of the space $\XQbfi$
to some domain $\oo\subset \Omega$, 
namely
\[
\norm{u}{L^2_{b,\phi}(0,\ell;L^2(\oo))}{} = \sup_{k\in\Indo} \phi^{1/2}(x_k)
\norm{u}{L^2(0,\ell;L^2(C_k))}{}.
\]
\par
In analogy with Theorems~\ref{teoequiv1}--\ref{teoequiv3}, one then proves:
\bete\label{teoequiv4}
Let $\phi$ be an admissible weight function of growth rate $\mu \in[0,1)$.
Then the function spaces $\XQbfi$, $\VQbfi$, $\ZQbfi$ and 
$\Lbfi{p}(0,\ell;L^p(\Omega))$ coincide
with the spaces $\tilde{L}^2_{b,\phi}(0,\ell;L^2(\Omega))$, 
$\tilde{L}^2_{b,\phi}(0,\ell;W^{1,2}(\Omega))$, 
$\tilde{L}^{2}_{b,\phi}(0,\ell;W^{-1,2}(\Omega))$, 
and $\tilde{L}^p_{b,\phi}(0,\ell;L^p(\Omega))$,
whose (equivalent) norms are given by
\begin{align}
\norm{u}{\tilde{L}^2_{b,\phi}(0,\ell;L^2(\Omega))}{2}=
&\sup_{\barx\in\Omega}\phi(\barx) \iq
    |u(x,t)|^2 \dwxt,
    \\
        \label{normVQ}
\norm{u}{\tilde L^{2}_{b,\phi}(0,\ell;W^{1,2}(\Omega))}{2}=
&\sup_{\barx\in\Omega}\phi(\barx) \iq
    \big( |u(x,t)|^2 + |\nabla u(x,t)|^2 \big) \dwxt,
    \\
        \label{dualQ}
\norm{u}{\tilde L^2_{b,\phi}(0,\ell;W^{-1,2}(\Omega))}{}=
& \sup_{v} \sup_{\barx\in\Omega}\phi(\barx) \iq
    u(x,t) v(x,t) \dwxt,
    \\
\norm{u}{\tilde{L}^p_{b,\phi}(0,\ell;L^p(\Omega))}{p}=
&\sup_{\barx\in\Omega}\phi(\barx) \iq
    |u(x,t)|^p \dwxt,
\end{align}
respectively. The supremum in \eqref{dualQ} is taken
over all $v$ such that the norm in 
$\tilde{L}^2_{b,\phi}(0,\ell;W^{1,2}(\Omega))$ 
is less or equal to $1$.
\ente
\begin{proof}
Omitted as being completely analogous to the three preceding
theorems. The only difference is an extra integration
over $t\in(0,\ell)$. Note that the already proven equivalences can
be simply written as
\begin{equation} \label{ekvxp}
\begin{aligned}
\norm{u}{\Lbfi{p}(0,\ell;L^p(\Omega))}{p}
&\approx \sup_{\barx} \phi(\barx) 
\norm{u}{L^p(0,\ell;L^p_{\barx}(\Omega))}{p},
\\
\norm{u}{\Lbfi{2}(0,\ell;W^{1,2}(\Omega))}{2}
&\approx \sup_{\barx} \phi(\barx) 
\norm{u}{L^2(0,\ell;W^{1,2}_{\barx}(\Omega))}{2},
\qquad \textrm{etc.}
\end{aligned}
\end{equation}

\end{proof}

\beos\label{notsame}
Note that in the above definitions, one \emph{first} integrates over
$t\in (0,\ell)$ and \emph{then} takes the weighted supremum. It is
thus clear that, e.g.,
\[
L^2(0,\ell; \ltwobfi ) \subset \XQbfi
\subset L^2_{\loc}(Q),
\]
where $Q=[0,\ell]\times \Omega$ and both inclusions are indeed strict. 
\eddos


\subsection{Some auxiliary results}

Given a precompact set $K$ in a metric space $M$,
we define Kolmogorov's $\varepsilon$-entropy as
\[
\hent{\varepsilon}{K}{M} := \ln N_{\varepsilon}(K,M),
\]
where $N_{\varepsilon}(K,M)$ is the smallest number of $\varepsilon$-balls
that cover $K$. Also, the symbol $B_r(u;M)$ denotes a ball centered in $u$,
of radius $r>0$, measured in the metric of $M$.

The following explicit version of the Aubin-Lions Lemma will be
instrumental in the proof of the main theorem.

\bele \label{le-aub_res}
Set
\begin{equation} \label{W-norm}
\norm{\chi}{\WQbfi}{} := \norm{\chi}{L^2_{b,\phi}(0,\ell;W^{1,2}(\Omega))}{} 
+ \norm{\dt\chi}{L^{2}_{b,\phi}(0,\ell;W^{-1,2}(\Omega))}{}.
\end{equation}

Let $\oo \subset \Omega$ be a ``reasonable'' domain in the sense
that
\begin{equation} \label{reason}
\# \Indo \le c_1 \vol(\oo).
\end{equation}
Let $r>0$, $\theta\in(0,1)$ be given. Then
\[
\hent{\theta r}{B_r(\chi;\WQbfi)}{L^2_{b,\phi}(0,\ell;L^2(\oo))} 
\le c_0 \vol(\oo);
\]
where the constant $c_0$ only depends on $c_1$, $\ell$ 
and $\theta$, but is independent
of $\chi$, $r$, $\oo$ and the weight function $\phi$ as long as
\eqref{weight-gr} and \eqref{weight-bla} are satisfied.
\enle
\begin{proof}
Observe that balls of radii $R\ge1$ are ``reasonable'' class
of domains and we will not work with any other $\oo$.

\noindent
STEP 1. Assume $\phi\equiv1$. Then $\WQbfi$ estimates 
from above each seminorm
\[
\norm{\chi}{L^2(0,\ell;W^{1,2}(C_k))}{}
+
\norm{\dt\chi}{L^2(0,\ell;W^{-1,2}(C_k))}{}, 
\qquad k \in \Indo.
\]
By the usual version of Aubin-Lions Lemma (see e.g. \cite{Si}),
we then have
\[
\hent{\theta r}{B_r(\chi;\WQbfi)}{L^2(0,\ell;L^2(C_k))}
\le c_1,
\]
where $c_1$ is independent of $k$. The desired covering arises as
a product of those and, in view of \eqref{reason}, the final estimate
follows.

\noindent
STEP 2. The case with general $\phi$ is reduced to the previous step
using the operator
\[
F:\chi \mapsto \phi^{1/2}\chi.
\]
The proof will be finished once we show that
\[
\norm{\chi}{N_{b,\phi}}{} \approx
\norm{F\chi}{N_{b,1}}{}
\]
and the equivalence constants can be taken independently 
on choosing $N_{b,\phi}$ as any of the spaces $L^2_{b,\phi}(0,\ell;L^2(\oo))$,
$\VQbfi$ or $\ZQbfi$.

(i) The case $N_{b,\phi}=L^2_{b,\phi}(0,\ell;L^2(\oo))$ clearly follows
from the fact that
\[
|F\chi(x,t)|^2 = \phi(x) |\chi(x,t)|^2 \approx \phi(x_k) |\chi(x,t)|^2,
\]
if $x\in C_k$.

(ii) Regarding the space $\VQbfi$, one obviously has
(cf.~\eqref{weight-bla})
\begin{equation} \label{ii-1}
|\nabla F\chi|^2 \le c_1 \phi \big( |\nabla \chi|^2 + |\chi|^2 \big)
\le c_2 \phi(x_k)  \big( |\nabla \chi|^2 + |\chi|^2 \big),
\end{equation}
for $x\in C_k$. The opposite inequality is more delicate. It is now
crucial that \eqref{weight-bla} holds with $1$, hence
\[
|\nabla F\chi| \ge |\phi^{1/2} \nabla \chi| - 
\frac12 \phi^{-1/2}|\nabla \phi| |\chi|
\ge \phi^{1/2}\big( |\nabla \chi| - \frac12 |\chi| \big)
.
\]
It then follows that
\[
|\nabla F\chi|^2 + |F\chi|^2 \approx
\big( |\nabla F\chi| + |F\chi| \big)^2
\ge 
c_3 \phi(x_k) \big( |\nabla \chi|^2 + |\chi|^2 \big)
\]
and the equivalence is concluded as in (i).

(iii) We first have to remark that in $\ZQbfi$ the operator $F$ is
defined by duality, i.e.,
\[
\dual{F\chi}{v} := \dual{\chi}{Fv}.
\]
But then
\begin{align*}
\norm{F\chi}{L^2_{b,1}(0,\ell;W^{-1,2}(\Omega))}{}
&= \sup_k \norm{F\chi}{L^2(0,\ell;W^{-1,2}(C_k))}{}
\\
&= \sup_k \sup_v \int_{C_k\times(0,\ell)} \chi \phi^{1/2}v \,dxdt
\\
&\approx \sup_k \sup_v \phi^{1/2}(x_k) \iq \chi v \,dx dt
\\
&
= \norm{\chi}{\ZQbfi}{}.
\end{align*}
Here the supremum is taken over all $v\in L^2(0,\ell;W^{1,2}_0(C_k))$
with unit norm; in the second step we have used the equivalence 
\begin{equation}
\norm{\phi^{1/2}v}{W^{1,2}(C_k)}{} \approx
\phi(x_k)^{1/2}
\norm{v}{W^{1,2}(C_k)}{}
\end{equation}
established in part (ii).

\end{proof}


\section{Well-posedness}
\label{secwell}

Here we give a rigorous mathematical formulation
of equation \eqref{eqref} within the spaces of 
uniformly locally $L^2$-functions. We first specify our basic
assumptions on the data, starting with the nonlinear diffusion
term:
\begin{align}\label{hpa1}
  & a\in C^0(\RR^d;\RR^d), \qquad
   a(0)=0,  \qquad
   \big(a(\xi)-a(\eta)\big)\cdot(\xi-\eta)
    \ge \kappa|\xi-\eta|^2, \quad\perogni\xi,\eta\in\RR^d,\\
 \label{hpa2}
  & \big|a(\xi)-a(\eta)\big| \le c\kappa|\xi-\eta|,
   \quad\perogni\xi,\eta\in\RR^d,\\
 \label{hpa3}
  &  \xi\mapsto a(\xi)\cdot \xi
   \quext{is a convex function on }\,\RR^d,
\end{align}  
where $\kappa>0$ and $c\ge 1$ are suitable constants.
We now introduce the family of nonlinear elliptic
operators $\{A_{\barx}\}_{\barx\in\RR^d}$ as
\begin{equation}\label{defiA}
  A_{\barx}: \wtwox \to W^{-1,2}_{\barx}(\Omega), \qquad
   \duav{A_{\barx}v,z}_{\barx} := 
     \io a(\nabla v(x)) \cdot 
       \Big( \nabla z(x) - z(x) \frac{x-\barx}{|x-\barx|} \Big)
       e^{-|x-\barx|}\,\dix.
\end{equation}  
In particular, if $v\in W^{1,2}_{b,\phi}(\Omega)$
for an admissible weight $\phi$ of growth rate $\mu < 1$, then 
$A_{\barx} v$ is an element of $W^{-1,2}_{\barx}(\Omega)$ for 
all $\barx\in\RR^d$.
The nonlinear function $f$ is assumed to satisfy
\begin{align}
\label{hpf0}
   f  \in C^0(\RR;\RR),& \quad  f(0)=0, 
\\ \label{hpf2}
  |f(r)-f(s)| & \le c_2(1 + |r|+|s|)^{p-2}|r-s|, \quad \perogni r,s\in \RR,
\\ \label{hpf3}
  \big( f(r) - f(s) \big)(r-s) &\ge -C |r-s|^2,\quad \perogni r,s\in \RR,
\\
 \label{hpf4}
   c_4 |r|^p - c_5 & \le f(r)r
   \le c_6 (|r|^p+1)\quad \perogni r\in \RR.
\end{align}  
for some $C$, $c_i>0$ and some $p\in(2,\infty)$.
Hence, we are requiring that $f$ grows superlinearly at
infinity, which holds in most applications.
As far as $h$ is concerned, we let
\begin{align}\label{hph}
  & \h:\Omega\times\RR^d\to\RR, \qquad
   \xi\mapsto \h(x,\xi) \quext{is globally
      Lipschitz for a.e.~}\,x\in\Omega,\\
  \label{hph2}
  & x\mapsto \h(x,\xi) \quext{is measurable and
    essentially bounded for all~}\,\xi\in\RR^d.
\end{align} 
\beos\label{alsoonu}
 With minor modifications in the proofs, 
 one could admit also a (Lipschitz) dependence
 on $u$ in the convective term $\h$.
 We limit ourselves to the slightly more restrictive
 setting \eqref{hph}-\eqref{hph2} just for the 
 sake of notational simplicity. 
\eddos
Finally, we take 
\begin{equation}\label{hpg}
  g \in \ltwob.
\end{equation}  
We are now  able to state our result on well-posedness
and dissipativity of the reaction-diffusion system
in the space $\ltwob$. %
Notice that, since we are only considering {\sl negative}\/
exponential weights, $L^2_b(\Omega)$ is continuously included 
into $\ltwobfi$ for any such weight. In particular, estimate
\eqref{diss} below makes sense.

\bete\label{teowell}
 Let assumptions\/ \eqref{hpa1}-\eqref{hpa3} and
 \eqref{hpf0}-\eqref{hpg} hold. Let also
 \begin{equation}\label{hpu0}
   u_0 \in \ltwob.
 \end{equation}  
 Then, there exists a\/ {\rm unique} function $u$ such that, for any
 $\barx \in \Omega$, one has
 \begin{equation}\label{regou}
 \begin{aligned}
   u &\in C^0([0,T];\ltwox) \cap L^2(0,T;\wtwox) 
   \cap L^{p}(0,T;\lpx),
    \\
   u_t &\in L^2(0,T;\wmtwox)+L^{p'}(0,T;\lppx), \quad
  \end{aligned}
 \end{equation}  
 and for all $\barx\in\Omega$ there holds
 \begin{equation}\label{eq}
   u_t + A_{\barx} u + f(u) + \h(\cdot,\nabla u) = g,
    \qquext{in }\,L^2(0,T;\wmtwox)+L^{p'}(0,T;L^{p'}_{\barx}(\Omega)).
 \end{equation}  
 Moreover, we have
 \begin{equation}\label{init}
   u|_{t=0}=u_0, \qquext{in }\, \ltwox.
 \end{equation}  
 Finally, for every admissible weight function
 $\phi$ with growth rate $\mu < 1$
 and almost all $t\ge 0$,
 there holds the dissipative estimate
 \begin{equation}\label{diss}
   \|u(t)\|_{\ltwobfi}^2
    + c_1 \norm{u}{\Lbfi{2}(t,t+1;W^{1,2}(\Omega))}{2}
    + c_2 \norm{u}{\Lbfi{p}(t,t+1;L^p(\Omega))}{p}
   \le \|u_0\|_{\ltwobfi}^2 e^{-\sigma t}
    + c_3,
 \end{equation}  
 where $\sigma$ and $c_i$ are positive constants depending
 on the parameters of the system, but independent of
 the initial datum $u_0$.
\ente
A function $u$ under the conditions of Theorem \ref{teowell} will
be simply called a ``solution'' in the sequel. Of course,
due to arbitrariness of $T$ any solution can be thought
to be defined for almost any $t\in(0,\infty)$.
\beos\label{whichspace}
 Equation \eqref{eq} can be also written in an
 expanded way as 
 \begin{align}\no
  & \io u_t(x,t) v(x,t) e^{-|x-\barx|}\,\dix
    + \io a(\nabla u(x,t)) \cdot
      \Big( \nabla v(x,t) - v(x,t) \frac{x-\barx}{|x-\barx|}\Big)
    e^{-|x-\barx|}\,\dix\\
  \label{eq2}
   & \mbox{}~~~~~~~~~~
   + \io \big(f(u(x,t))+\h(x,\nabla u(x,t))
            -g(x)\big) v(x,t) e^{-|x-\barx|}\,\dix
    = 0,
 \end{align}
 the above being intended to hold for any
 $\barx\in\Omega$, almost any $t\in(0,T)$ and 
 any test function 
 $v \in L^2(0,T;\wtwox)\cap L^{p}(0,T;L^{p}_{\barx}(\Omega))$.
 In particular, by \eqref{regou}, one can take $v=u$.
\eddos
\begin{proof}
The proof is carried out by suitably
approximating \eqref{eq} through a family of
problems defined on bounded
domains and then passing to the limit via monotonicity
and compactness methods.

As a first step, we then define $\Omega_n:=B_n(0,\RR^d)$, 
$n\in\NN$, and for any 
$n$ consider a cutoff function 
$\psi_n\in C^\infty(\Omega;[0,1])$ such that $\psi\equiv 1$
in $\bar \Omega_{n-1}$ and $\supp(\psi)\subset \Omega_n$.
Then, we set $u\zzn:= u_0\psi_n$ and $g_n:=g\psi_n$.
Thanks to \eqref{hpg} and \eqref{hpu0},
applying Lebesgue's theorem
one can easily check that, for every $\barx\in\Omega$,
\begin{equation}\label{cou0g}
  u\zzn \to u_0 \quand g_n\to g,
   \quext{strongly in }\,\ltwox.
\end{equation}  
We also set $X_n:=L^2(\Omega_n)$,
$V_n:=W^{1,2}_0(\Omega_n)$ and
define the elliptic operator
\begin{equation}\label{defiA2}
  A_n: V_n\to V_n', \qquad
   \duav{A_n v,z}:=\io \nabla a(v(x))\cdot \nabla z(x)\,\dix,
\end{equation}  
where $v,z\in V_n$.
Then, we can introduce our approximate problem
\begin{align}\label{eqn}
  & u_{n,t} + A_n u_n + f(u_n) + \h(\cdot,\nabla u_n) = g_n,
   \quext{in }\,L^2(0,T;V_n)+L^{p'}(0,T;L^{p'}(\Omega_n)),\\
  \label{inizn}
  & u_n|_{t=0}=u\zzn,
   \quext{a.e.~in }\,\Omega_n.
\end{align} 
We have the following
\bele
 For all $n\in\NN$, there exists one and only one solution 
 $u_n$ to\/ \eqref{eqn}-\eqref{inizn} such that
 \begin{equation}\label{regoun}
   u_{n,t} \in L^2(0,T;V_n')+L^{p'}(0,T;L^{p'}(\Omega_n)), \quad
   u_n \in C^0([0,T];X_n) \cap L^2(0,T;V_n) \cap L^{p}(0,T;L^p(\Omega_n)).
 \end{equation}  
\enle
\noindent%
The proof of the lemma is more or less standard and
mainly relies on the basic tools of the theory of maximal monotone
operators. We will not give it since most of the difficulties
will be the same we will face in the passage to the limit
$n\nearrow\infty$ we now describe.

Assume $u_n$ be extended to $0$ outside $\Omega_n$ 
and test \eqref{eqn} by $u_n e^{-|\cdot-\barx|}$,
for arbitrary $\barx\in \Omega$. Then, we readily obtain the
basic estimate
\begin{align}\no
  & \frac12\ddt\io |u_n(x,t)|^2 e^{-|x-\barx|}\,\dix
   + \io a(\nabla u_n(x,t)) \cdot \Big(
      \nabla u_n(x,t) - u_n(x,t) \frac{x-\barx}{|x-\barx|} \Big)
      e^{-|x-\barx|}\,\dix\\
  \label{st}
  & \mbox{}~~~~~~~~~~  
    + \io \big( f(u_n(x,t)) + \h(x,\nabla u_n(x,t)) \big) 
         u_n(x,t) e^{-|x-\barx|}\,\dix
    = \io g_n(x) u_n(x,t) e^{-|x-\barx|}\,\dix.
\end{align}
Using hypotheses \eqref{hpa1}-\eqref{hpa2} and 
\eqref{hpf2}, it is then 
not difficult to deduce from
\eqref{st} a priori estimates
in weighted spaces which entail
\begin{equation}\label{co11}
  u_n\to u \quext{weakly star in }\,
   L^\infty(0,T;\ltwox)\cap L^p(0,T;L^p_{\barx}(\Omega))
    \cap L^2(0,T;\wtwox).
\end{equation}
Note that, here and below, all convergence relations
are intended up to the extraction of subsequences,
not relabelled (see also Remark~\ref{qualebarx} below
for more details).
Next, writing \eqref{eqn} in the form corresponding
to \eqref{eq2}, namely
\begin{align}\no
 & \io u_{n,t}(x,t) v(x,t) e^{-|x-\barx|}\,\dix
   + \io a(\nabla u_n(x,t)) 
     \cdot \Big( \nabla v(x,t) - v(x,t) \frac{x-\barx}{|x-\barx|}\Big)
   e^{-|x-\barx|}\,\dix\\
 \label{eq2n}
  & \mbox{}~~~~~~~~~~
  + \io \big(f(u_n(x,t))+\h(x,\nabla u_n(x,t))-g_n(x)\big)
         v(x,t) e^{-|x-\barx|}\,\dix
   = 0,
\end{align}
and letting $v$ vary in 
$L^2(0,T;\wtwox) \cap L^{p}(0,T;L^p_{\barx}(\Omega))$,
passing to the supremum with respect to $v$ of unit norm,
it is not difficult to obtain
\begin{equation}\label{co12}
  u_{n,t}\to u_t \quext{weakly in }\,
   L^2(0,T;\wmtwox) + L^{p'}(0,T;L^{p'}_{\barx}(\Omega)).
\end{equation}
At this point, if one considers the restrictions to a fixed
domain $\Omega_m$, then \eqref{co11} implies in 
particular
\begin{equation}\label{co11m}
  u_n\to u \quext{weakly star in }\,
   L^\infty(0,T;X_m)\cap L^p(0,T;L^p(\Omega_m))
    \cap L^2(0,T;W^{1,2}(\Omega_m)).
\end{equation}
On the other hand, if we write \eqref{eqn} for 
$n>m$ and test it by a generic
$v\in L^2(0,T;V_m)\cap L^p(0,T;L^p(\Omega_m))$
(extended by $0$ outside $\Omega_m$), then,
using \eqref{co11m} and applying duality arguments,
we readily infer
\begin{equation}\label{co12m}
  u_{n,t}\to u_t \quext{weakly in }\,
   L^2(0,T;V_m') + L^{p'}(0,T;L^{p'}(\Omega_m)).
\end{equation}
In particular, by the Aubin-Lions Lemma, 
we get from \eqref{co11m}-\eqref{co12m} that
\begin{equation}\label{co13m}
  u_n\to u \quext{strongly in }\,
   L^2(0,T;X_m).
\end{equation}
More precisely, by arbitrariness of 
$m$, we have 
\begin{equation}\label{co13}
  u_n\to u \quext{a.e.~in }\,
   \Omega\times (0,T).
\end{equation}
Thus, recalling \eqref{co11} and
applying Lebesgue's Theorem with respect to the 
measure $\deriv_{\barx}\!x=e^{-|x-\barx|}\deriv\!x$
(notice that $\Omega=\RR^d$ has {\sl finite}\/
$\deriv_{\barx}\!x$-measure),
we readily obtain
\begin{equation}\label{co14}
  u_n\to u \quext{strongly in }\/
   L^q(0,T;L^q_{\barx}(\Omega)), \quad\perogni
   q\in [1,p)
\end{equation}
and, thanks to \eqref{hpf2},
\begin{equation}\label{co15}
  f(u_n)\to f(u) \quext{strongly in }\/
   L^q(0,T;L^q_{\barx}(\Omega)), \quad\perogni
   q\in [1,p').
\end{equation}
%
%
Thus, we are now ready to pass to the limit in equation
\eqref{eqn}. To do this, we first observe that, 
by \eqref{co11} and assumptions \eqref{hpa2}
and \eqref{hph}, if $\barx\in\Omega$
is fixed, there exist $\aalpha\in L^2(0,T;\ltwox^d)$
and $\boldsymbol{\h}\in L^2(0,T;\ltwox)$
such that 
\begin{align}\label{co16}
  & a(\nabla u_n)\to \aalpha \quext{weakly in }\/
   L^2(0,T;\ltwox^d),\\
 \label{co16h}
  &  \h(\cdot,\nabla u_n)\to \tilde{\h} \quext{weakly in }\/
   L^2(0,T;\ltwox). 
\end{align}
Notice that, a priori, $\aalpha$ 
and $\tilde{\h}$ might depend on the 
choice of $\barx$. Let us now come back to \eqref{eq2n}.
It is clear that we can take its limit, which assumes
the form
\begin{align}\no
  & \io u_{t}(x,t) v(x,t) e^{-|x-\barx|}\,\dix
   + \io \aalpha(x,t)\cdot
     \Big( \nabla v(x,t) - v(x,t) \frac{x-\barx}{|x-\barx|}\Big)
   e^{-|x-\barx|}\,\dix\\
 \label{eq2lim}
  & \mbox{}~~~~~~~~~~
  + \io \big(f(u(x,t))+\tilde{\h}(x,t)-g(x)\big) 
    v(x,t) e^{-|x-\barx|}\,\dix
   = 0.
\end{align}
Now, let us choose $v=u_n$ in \eqref{eq2n}, rearrange some terms,
integrate over $(0,T)$, and take the supremum limit. This 
procedure gives
\begin{align}\no
  & \limsup_{n\nearrow\infty} \iTo a(\nabla u_n(x,t))
   \cdot  \nabla u_n(x,t) e^{-|x-\barx|}\,\dix\,\dit\\
 \no
  &\mbox{}~~~~~ 
   \le -\frac12 \liminf_{n\nearrow\infty} \io |u_n(x,T)|^2 e^{-|x-\barx|}\,\dix
   + \frac12 \limsup_{n\nearrow\infty} \io |u\zzn(x)|^2 e^{-|x-\barx|}\,\dix\\
 \no 
  & \mbox{}~~~~~
   - \liminf_{n\nearrow\infty} \iTo \big(f(u_n(x,t))+C u_n(x,t)\big) 
      u_n(x,t) e^{-|x-\barx|}\,\dix\,\dit\\
 \no 
  & \mbox{}~~~~~
   + \limsup_{n\nearrow\infty} \iTo C |u_n(x,t)|^2 e^{-|x-\barx|}\,\dix\,\dit
   + \limsup_{n\nearrow\infty} \iTo g_n(x) u_n(x,t) e^{-|x-\barx|}\,\dix\,\dit\\
 \no
  & \mbox{}~~~~~ 
  - \liminf_{n\nearrow\infty} \iTo \h(x,\nabla u_n(x,t)) 
   u_n(x,t) e^{-|x-\barx|}\,\dix\,\dit\\
 \label{giroeq}
  & \mbox{}~~~~~ 
  + \limsup_{n\nearrow\infty} \iTo a(\nabla u_n(x,t)) \cdot
     \frac{x-\barx}{|x-\barx|} u_n(x,t) e^{-|x-\barx|}\,\dix\,\dit.
\end{align}
At this point, we aim to compute the limits on the \rhs.
First, let us observe that the first two terms are treated
by means of \eqref{cou0g}, \eqref{co11}, and semicontinuity of
norms with respect to weak star convergence.
Next, recalling \eqref{hpf0} and \eqref{hpf3}, by \eqref{co13}
and Fatou's Lemma we obtain
\begin{align}\no
  & \iTo \big(f(u(x,t))+C u(x,t)\big) u(x,t) 
     e^{-|x-\barx|}\,\dix\,\dit \\
 \label{co17}
  & \mbox{}~~~~~
   \le \liminf_{n\nearrow\infty} 
    \iTo \big(f(u_n(x,t))+C u_n(x,t)\big) 
       u_n(x,t) e^{-|x-\barx|}\,\dix\,\dit.
 \end{align}
The subsequent three terms are treated thanks to \eqref{cou0g},
\eqref{co14} (where we can take $q=2$)
and \eqref{co16h}. Finally,
using \eqref{co16} and again \eqref{co14}, we arrive at
\begin{align}\no
  & \lim_{n\nearrow\infty} \iTo a(\nabla u_n(x,t))\cdot
    \frac{x-\barx}{|x-\barx|} u_n(x,t) e^{-|x-\barx|}\,\dix\,\dit\\
 \label{co18}
  & \mbox{}~~~~~
    = \iTo \aalpha(x,t) \cdot \frac{x-\barx}{|x-\barx|} 
        u(x,t) e^{-|x-\barx|}\,\dix\,\dit.
\end{align}
Thus, comparing \eqref{giroeq} with \eqref{eq2lim} 
(written for $v=u$ and integrated in time), 
we finally deduce that
\begin{align}\no
  & \limsup_{n\nearrow\infty} \iTo a(\nabla u_n(x,t))
   \cdot  \nabla u_n(x,t) e^{-|x-\barx|}\,\dix\,\dit\\
 \label{co19}
  & \mbox{}~~~~~
    \le \iTo \aalpha(x,t) \cdot \nabla u(x,t) e^{-|x-\barx|}\,\dix\,\dit.
\end{align}
Noting now that, by assumption \eqref{hpa1}, 
$a$ induces a maximal monotone operator on the 
Hilbert space $L^2(0,T;\ltwox^d)$, the usual
monotonicity argument (cf., e.g., 
\cite[Prop.~1.1, p.~42]{Ba}) 
permits to say that
\begin{equation}\label{co20}
  \aalpha(x,t)=a(\nabla u(x,t))
   \quext{$\deriv_{\barx}\!x$-a.e.~in }\,\Omega
   ~~\text{and a.e.~in }\,(0,T),
\end{equation}
whence the same holds 
almost everywhere with respect to Lebesgue's
measure in $\Omega\times (0,T)$. In particular, 
$\aalpha$ is independent of the choice of 
$\barx$. Thus, substituting in \eqref{eq2lim},
we get exactly \eqref{eq2}.
Finally, we notice that, as a consequence of
\eqref{co19}-\eqref{co20} and lower semicontinuity,
\begin{equation} \label{co19bis}
  \iTo a(\nabla u_n(x,t)) \cdot  \nabla u_n(x,t) e^{-|x-\barx|}\,\dix\,\dit
   \to \iTo a(\nabla u(x,t)) \cdot \nabla u(x,t) e^{-|x-\barx|}\,\dix\,\dit.
\end{equation}
Thus, using \eqref{hpa3} and, e.g., \cite[Thm.~2.11]{Fe},
we obtain
\begin{equation} \label{co19.3}
  \nabla u_n(x,t) \to \nabla u(x,t) 
   \quext{a.e.~in }\,\Omega\times (0,T),
\end{equation}
whence, by \eqref{hph}, \eqref{co16} and Lebesgue's
Theorem,
\begin{equation} \label{co19.4}
  \nabla u_n \to \nabla u~~\text{and }\,
  \h(\cdot,\nabla u_n) \to \h(\cdot,\nabla u)
   \quext{strongly~in }L^q(0,T;L^q_{\barx}(\Omega))
\end{equation}
for all $q\in[1,2)$. In particular, 
$\tilde{\h}=\h(\cdot,\nabla u)$ (cf.~\eqref{co16h}),
which concludes the proof of existence.
\beos\label{qualebarx}
 It is worth observing that relations 
 \eqref{co11}-\eqref{co12} and \eqref{co14}-\eqref{co15}
 hold for {\sl any}\/ $\barx\in\RR^d$ and the limits are 
 independent of $\barx$. 
 This follows already from the
 fact that the spaces $\lpx$ coincide for different
 values of $\barx$. However, it is still necessary to
 consider the weak formulation for all $\barx$ simultaneously 
 to make sure that the a priori estimates are also uniform
with respect to $\barx$. In virtue of the equivalence relations
 \eqref{ekvx} and \eqref{ekvxp} this then leads
 to the estimates in the uniformly bounded spaces.
\eddos
\beos\label{convexity}
 In the case when $\h$ is a {\sl linear}\/ convection
 term (namely, $\h(x,\xi)={\boldsymbol v(x)}\cdot \xi$
 for some measurable and bounded function ${\boldsymbol v}$),
 then assumption \eqref{hpa3} can be avoided. Actually,
 the only role of \eqref{hpa3} is that of guaranteeing
 the strong convergence \eqref{co19.4} of gradients,
 which is not required for taking the limit in case
 $\h$ is linear.
\eddos
\beos\label{vector1}
It is not difficult to realize that Theorem~\ref{teowell} can be extended to systems of $m$ equations
provided that the nonlinear function $a$ is replaced by $\mathbf{a}\in C^0(\RR^{m\times d};\RR^{m\times d})$
satisfying suitable reformulations of \eqref{hpa1} and \eqref{hpa2} and $h$ is replaced by a linear function
of the form $\mathbf{h}(x,{\mathbf M})=\boldsymbol{v}(x)\cdot {\mathbf M}$, where $\mathbf{M} \in \RR^{m\times d}$ (see Remark~\ref{convexity}).
Another possibility is to preserve a nonlinear convective term $\mathbf{h}:\Omega\times\RR^{m\times d}\to\RR^m$ 
satisfying suitable generalizations of \eqref{hph} and \eqref{hph2} and taking the vector Laplacian
$-\mathbf{\Delta}$ as diffusion operator.
\eddos
\medskip

Let us now move to dissipativity. To prove it,
let us go back to \eqref{eq2lim}, take $v=u$ 
and use \eqref{co20}, \eqref{hpa1}-\eqref{hpa2},
\eqref{hph} and \eqref{hpf2}. Then, we deduce, for some 
$\sigma>0$ independent of $\barx$,
\begin{equation}\label{co31}
  \ddt \|u\|_{\ltwox}^2
   + \sigma \big(
   \| \nabla u \|_{\ltwox}^2
   + \| u \|_{L^p_{\barx}(\Omega)}^p \big)
   \le c \io \big(1+g^2(x)\big)e^{-|x-\barx|}\,\dix
   \le c_1,
\end{equation}
where $c_1$ only depends on $\norm{g}{\ltwob}{}$
and on the Lipschitz constant of $\h$
(and is independent of $\barx$). By a standard application
of Gronwall's lemma, we further deduce
\[
\norm{u(t)}{\ltwox}{2}
+ \sigma \int_t^{t+1} \big(
\norm{\nabla u}{\wtwox}{2}
+ \norm{u}{\lpx}{p} \big) \dis
\le 
\norm{u_0}{\ltwox}{2} e^{-\sigma t} + c_2
.
\]
Next, we multiply with $\phi(\barx)$, and take the supremum
over $\barx \in \RR^d$. Using the fact that $\phi$ is uniformly
bounded and also the equivalence relations \eqref{ekvx}, 
\eqref{ekvxp}, we finally conclude \eqref{diss}.

\par
Notice that the above is a dissipative
estimate in any of the spaces $\ltwobfi$
where $\phi$ is an admissible weight of growth rate
strictly lower than $1$.

\medskip

Finally, let us prove uniqueness, which is standard. 
Indeed, it is sufficient to write \eqref{eq} for
a couple of solutions $u_1$ and $u_2$,
take the difference, test it by $u_1-u_2$
(in the appropriate functional sense)
and integrate with respect to the measure 
$\deriv_{\barx}\!x\otimes\dit$. Then the thesis follows
as before by using Gronwall's lemma and taking
the supremum  with respect to $\barx$. We omit the details
since we shall prove more refined contractive estimates
in the next section (Theorem~\ref{te-cont1}).
The proof of Theorem~\ref{teowell} is complete.
\end{proof}
In order to prepare the long time analysis, we need a further regularity
result.
\bete\label{th:reg} 
 Let $d\leq 3$ and consider a solution $u$. 
 Then, for any $q\in(1,\infty)$ and any $\tau>0$,
 $u$ enjoys the additional regularity
 \begin{equation} \label{reg2}
   u \in L^{\infty}(\tau,\infty;L^{q}_b(\Omega)).
 \end{equation}
 More precisely, for any $q\in(1,\infty)$ 
 there exists a computable 
 nonnegative-valued function $\mathcal{Q}$,
 depending on $q$ and increasingly
 monotone in each of its arguments, such that
 \begin{equation} \label{reg2.2}
   \|u(t)\|_{L^{q}_b(\Omega)}
    \le \mathcal{Q}(\tau^{-1},\|u_0\|_{L^2_b(\Omega)}),
    \qquad\perogni t\ge\tau>0.
 \end{equation}
\ente
\begin{proof}
The proof is performed by means of (finitely many) 
iterative estimates. As a first step, we notice that,
due to the dissipative estimate \eqref{diss},
for any $t\ge 0$, any $\tau\in(0,1)$ 
and any $\barx\in\Omega$ there exists 
$t_0\in[t,t+\tau]$ (possibly depending also 
on $\barx$) such that
\begin{equation}\label{consdiss1}
  \|u(t_0)\|_{L^{p}_{\barx}(\Omega)}
   +\|u\|_{L^2_{b}(t,t+1;W^{1,2}(\Omega))}
    \le \mathcal{Q}(\tau^{-1},\|u_0\|_{L^2_b(\Omega)}),
\end{equation}
where $\mathcal{Q}$ is as in the statement.

Then, we can test the equation by $v=|u|^\alpha u$, 
where $\alpha=p-2>0$ due to our
assumptions. Such a test function is indeed admissible at least
on the level of approximations, thanks to uniqueness.
Thus, using \eqref{hpf4} and \eqref{hph}, and
observing that $a(\nabla u) \cdot \nabla v \ge 0$
by \eqref{hpa1}, one deduces after obvious manipulations
\begin{equation} \label{conto41}
  \ddt \frac1p \norm{u}{L^p_{\barx}(\Omega)}{p}
   + c_1 \norm{u}{L^{p+\alpha}_{\barx}(\Omega)}{p+\alpha}
  \le c_2 + c_3 \norm{u}{L^p_{\barx}(\Omega)}{p}
  + c_4 \io \big(1+|\nabla u|+|g|\big) |u|^{\alpha+1} \dwx.
\end{equation}
The last integrand in \eqref{conto41}
is then simply estimated as
\[
  \big(|\nabla u|+|g|\big)|u|^{a+1} 
   \le \varepsilon |u|^{2\alpha+2} 
   + \varepsilon^{-1} \big( |\nabla u|^2 + |g|^2 \big).
\]
Then, choosing $\varepsilon$ small enough and 
remarking that $2a+2 = p+a$, we further deduce that
\begin{equation} \label{conto42}
  \ddt  \norm{u}{L^p_{\barx}(\Omega)}{p}
  + \frac{c_1}2 \norm{u}{L^{p+\alpha}_{\barx}(\Omega)}{p+\alpha}
  \le c_5 + c_6  \norm{u}{L^p_{\barx}(\Omega)}{p}
   + c_7 \norm{\nabla u}{L^2_{\barx}(\Omega)}{2}.
\end{equation}
Then, 
we can integrate \eqref{conto42}
over $(t_0,t_0+2)$. Recalling \eqref{consdiss1} and
using Gronwall's Lemma, we can then pass to the 
supremum with respect to $\barx$ first on the \rhs\ and
then on the \lhs. Noting that for any $\barx$ it
is $t_0(\barx)\le t+\tau$, by arbitrariness of $t$
in $\RR^+$ we deduce
\begin{equation}\label{consdiss2}
  \|u\|_{L^\infty(\tau,\infty;L^{p}_b(\Omega))}
    \le \mathcal{Q}(\tau^{-1},\|u_0\|_{L^2_b(\Omega)}).
\end{equation}
Moreover, we also obtain that, for each $t\ge \tau$
and any $\barx\in\Omega$,
there exists $t_1\in [t,t+\tau]$ such that
\begin{equation}\label{consdiss3x}
  \|u(t_1)\|_{L^{p+\alpha}_{\barx}(\Omega)}
   +\|u\|_{L^{p+\alpha}_{b}(t,t+1;L^{p+\alpha}(\Omega))}
    \le \mathcal{Q}(\tau^{-1},\|u_0\|_{L^2_b(\Omega)}).
\end{equation}

\smallskip

We can now proceed by an induction argument. 
More precisely, we will just need a finite number of 
steps. Actually, since $\alpha=p-2>1$, we will
stop after $n$ iterations when $n\in \NN$ is such that
$p+(n-1)\alpha=n\alpha+2\ge q$. 

So, we can assume that, given $k\le n$,
for each $t\ge \tau$ and any $\barx\in\Omega$,
there exists $t_{k-1}\in [t,t+\tau]$ such that
\begin{equation}\label{consdiss3}
  \|u(t_{k-1})\|_{L^{(k-1)\alpha+2}_{\barx}(\Omega)}
   +\|u\|_{L^{k\alpha+2}_{b}(t,t+1;L^{k\alpha+2}(\Omega))}
    \le \mathcal{Q}(\tau^{-1},\|u_0\|_{L^2_b(\Omega)}),
\end{equation}
and prove now the same relation with $k-1$ replaced
by $k$. 

To do this, we test the equation by $v=|u|^{k\alpha} u$, 
where $\alpha=p-2>0$ as before. Then, we obtain the analogue
of \eqref{conto41}, where, however, we need to use 
\eqref{hpa1} a bit more precisely. Namely, we get
\begin{align} \no
  & \ddt c_{1,k} \norm{u}{L^{k\alpha+2}_{\barx}(\Omega)}{k\alpha+2}
   + c_{2,k} \io |\nabla u|^2 |u|^{k\alpha} \dwx
   + c_{3,k} \norm{u}{L^{p+k\alpha}_{\barx}(\Omega)}{p+k\alpha}\\
 \label{conto51} 
  & \mbox{}~~~~~
  \le c_{4,k} + c_{5,k} \norm{u}{L^{k\alpha+2}_{\barx}(\Omega)}{k\alpha+2}
  + c_{6,k} \io \big(1+|\nabla u|+|g|\big) |u|^{k\alpha+1} \dwx.
\end{align}
All constants $c$ or $c_{i,k}$ here and below will be allowed
to depend on $k$. However, since a finite number of induction steps
will suffice, we will not need to compute them explicitly.
To estimate the terms on the \rhs, we then observe that
\begin{equation} \label{conto52}
  c_{6,k} \io \big(1+|\nabla u|\big) |u|^{k\alpha+1} \dwx
   \le \epsilon \io |\nabla u|^2 |u|^{k\alpha} \dwx
   + c_\epsilon
   + c_\epsilon \norm{u}{L^{k\alpha+2}_{\barx}(\Omega)}{k\alpha+2}.
\end{equation}
As for the $g$-term, we need however to be much more accurate
than before. Firstly, we notice that, 
for positive $\lambda_i$, $i=1,2,3$,
such that $\lambda_1+\lambda_2+\lambda_3=1$
(and that will be chosen below), we have
\begin{align} \no
  & c_{6,k} \io |g| |u|^{k\alpha+1} \dwx
   = c_{6,k} \io \big( |g| e^{-\lambda_1|x-\barx|} \big)
    \big( |u|^{\frac{3k\alpha}4} e^{-\lambda_2|x-\barx|} \big)
    \big( |u|^{\frac{k\alpha+4}4} e^{-\lambda_3|x-\barx|} \big) \,\dix \\
 \no
   & \mbox{}~~~~~~~\le c \big\| |g| e^{-\lambda_1|\cdot-\barx|} \big\|_{L^2(\Omega)} \times
   \big\| |u|^{\frac{3k\alpha}4} e^{-\lambda_2|\cdot-\barx|} \big\|_{L^{q_2}(\Omega)} \times
   \big\| |u|^{\frac{k\alpha+4}4} e^{-\lambda_3|\cdot-\barx|} \big\|_{L^{q_3}(\Omega)}\\[1mm]
 \label{conto53}  
   & \mbox{}~~~~~~~ =: \calI_1 \times \calI_2 \times \calI_3,
    \qquext{~~~where~}\,\frac{1}{q_2}+\frac{1}{q_3}=\frac12.
\end{align}   
Let us now estimate the quantities $\calI_i$. Actually, taking
\begin{equation}\label{sceltapq}
  q_2=\frac{4(k\alpha+2)}{k\alpha}, \qquad
   q_3=\frac{4(k\alpha+2)}{k\alpha+4},
\end{equation}
it is not difficult to obtain
\begin{equation}\label{stimaX}
  \calI_2=\Big\| |u|^{\frac{k\alpha+2}{2}} e^{-\frac{2(k\alpha+2)\lambda_2|\cdot-\barx|}{3k\alpha}} 
      \Big\|_{L^6(\Omega)}^{\frac{3k\alpha}{2(k\alpha+2)}},
\end{equation}
whence, by continuity of the embedding $H^1(\Omega)\subset L^6(\Omega)$,
it is straightforward  to arrive at 
\begin{equation}\label{stimaX2}
  \calI_2
   \le 
   c  
   \bigg| \io |u|^{k\alpha} |\nabla u|^2 e^{-\frac{4(k\alpha+2)\lambda_2|x-\barx|}{3k\alpha}}
    \bigg|^{\frac{3k\alpha}{4(k\alpha+2)}}  
   +  
   c 
   \bigg| \io |u|^{k\alpha+2} e^{-\frac{4(k\alpha+2)\lambda_2|x-\barx|}{3k\alpha}}
    \bigg|^{\frac{3k\alpha}{4(k\alpha+2)}}.
\end{equation}
Computing $\calI_3$ directly, we similarly obtain
\begin{equation}\label{stimaY}
  \calI_3 \le \bigg| \io |u|^{k\alpha+2} 
          e^{-\frac{4(k\alpha+2)\lambda_3|x-\barx|}{k\alpha+4}} 
    \bigg|^{\frac{k\alpha+4}{4(k\alpha+2)}}.
\end{equation}
At this point, in order to get the same weight functions as on the \lhs,
we choose 
\begin{equation}\label{sceltala}
  \lambda_2=\frac{3k\alpha}{4(k\alpha+2)}, \quad 
   \lambda_3=\frac{k\alpha+4}{4(k\alpha+2)}, 
   \qquext{so that }\,\lambda_1=\frac{1}{k\alpha+2}
\end{equation}
and consequently we obtain  
\begin{equation}\label{stimaI1}      
  \calI_1 \le c \bigg|\io g^2(x) 
       e^{-\frac{2|x-\barx|}{k\alpha+2}} \,\dix\bigg|^{1/2}
   \le c \sup_{\barx\in\Omega}\bigg\{
     \bigg|\io g^2(x) 
       e^{-\frac{2|x-\barx|}{k\alpha+2}} \,\dix\bigg|^{1/2}
        \bigg\} \le c,
\end{equation}
where the last inequality follows from the fact that 
we have obtained a norm of $g$ that is equivalent to the 
usual norm of $L^2_b(\Omega)$
(this fact can be shown proceeding similarly with the
proof of Theorem~\ref{teoequiv1} in
the case $\phi\equiv1$). Notice that, the larger is 
$k$, the slower is the decay of the exponential weight
(however, we will not need to take $k\to \infty$ here).
  
Thus, using also the Young inequality, 
\eqref{conto53} gives  
\begin{equation}\label{fineconto}
  \calI_1 \times \calI_2 \times \calI_3 
   \le \epsilon \bigg(
     \io |u|^{k\alpha} |\nabla u|^2 e^{-|x-\barx|}
     + \io |u|^{k\alpha+2} e^{-|x-\barx|} \bigg)
    + c_\epsilon  
     \bigg| \io |u|^{k\alpha+2} e^{-|x-\barx|} \bigg|^{\frac{k\alpha+4}{k\alpha+8}},
\end{equation}     
where of course the latter exponent is (strictly)
lower than $1$.

Thus, integrating \eqref{conto51} over $(t_{k-1},t_{k-1}+2)$,
taking $\epsilon$ small enough,
using Gronwall's Lemma, and taking as before the supremum
with respect to $\barx$ first on the \rhs\ and then on the \lhs,
it is almost immediate to obtain \eqref{consdiss3}
with $k-1$ replaced by $k$. This concludes the proof.
\end{proof}
\beos\label{dim}

Note that should we assume $g\in L^{\infty}(\Omega)$,
the Theorem~\ref{th:reg}  can be proved in a simpler way and the
restriction $d\le 3$ can be removed.
\eddos 
\beos\label{vector2}
It is not difficult to realize that Theorem~\ref{th:reg} can be extended to systems of $m$ equations
provided that the diffusion operator is the vector Laplacian
$-\mathbf{\Delta}$ and the nonlinear convective term $h$ is replaced by 
$\mathbf{h}:\Omega\times\RR^{m\times d}\to\RR^m$ 
satisfying suitable generalizations of \eqref{hph} and \eqref{hph2} (see also Remark~\ref{vector1}).
\eddos


\section{Global attractor}
\label{secdynam}

Thanks to Theorem~\ref{teowell} we can introduce the solution
operator
\begin{equation*}
  S(t):\ltwob \to \ltwob,\qquad u_0  \mapsto u(\cdot,t).
\end{equation*}
Before showing that $S(\cdot)$ 
is a continuous semigroup, we prove a simple
\beco\label{new:cor}
 The semiflow $S(\cdot)$ admits 
 an absorbing set of the form
 \begin{equation} \label{abs0}
   \bnul := B_{K}(0;\ltwob),
 \end{equation}
 with a sufficiently large $K>0$. Moreover, 
 $\bnul$ can be chosen to be positively invariant  
 and bounded in the space $L^{q}_b(\Omega)$
 for $q$ arbitrarily large.
\enco
\begin{proof}
The existence of an absorbing set $\bnul_0$ 
satisfying \eqref{abs0} is an immediate consequence
of \eqref{diss}. Setting
\begin{equation} \label{absnew}
  \bnul := \bigcup_{t\ge 1}S(t)\bnul_0,
\end{equation}
we immediately obtain the positive invariance, 
as well as the $L^{q}_b$-boundedness,
thanks to \eqref{reg2}.
\end{proof}
Notice that, however, we cannot expect that the
dynamics be compact in $\ltwob$. The standard way out of this impasse
is the local topology $\ltwoloc$. Indeed,
thanks to Corollary~\ref{new:cor}, we can restrict
our analysis to those trajectories taking values, for all
nonnegative times in the $L^{q}_b$-bounded set $\bnul$.
Then, one easily verifies that
\begin{equation} \label{loc-x}
 \begin{aligned}
  u_n \to u_0 \quad \textrm{in $\ltwoloc$}
   \qquad &\iff \qquad
  u_n \to u_0 \quad \textrm{in $\ltwox$},\\
  \end{aligned}
\end{equation}
for any $u_n$, $u_0\in \bnul$.
Namely, whatever is $\barx\in\Omega$, the norm
of $\ltwox$ induces to $\bnul$ exactly the 
$\ltwoloc$-topology (in particular, one
could directly choose $\barx=0$ at this stage).
Thus, recalling also that the 
solutions are continuous as functions
with values in $\ltwox$, this space seems to be most convenient
for the construction of the global attractor. More precisely, we
are going to establish the existence of the $(\ltwob,\ltwoloc)$-attractor,
following the terminology of \cite{BaVi92}.

We recall that one possible strategy to show the compactness of the dynamics in $\ltwoloc$
is to derive higher regularity estimates, as for example in $\wtwob$.
However, as we mentioned in the Introduction, 
here we adopt a more elementary approach, which circumvents
more advanced regularity techniques, resting only on the natural
parabolic compactness of solutions. This is easy to obtain while
we look at the dynamics from the perspective of ``trajectories''
with some finite fixed length $\ell$.

We then introduce the set of the {\sl short trajectories}\/
taking values in $\bnul$:
\[
  \xxQ := \big\{ \chi \in L^2(0,\ell;\ltwox);\ 
   \textrm{$\chi$ is a solution of  \eqref{eqref}, $\chi(0) \in \bnul$} \big\}.
\]
Further, we define the semigroup
\begin{equation*}
 L(t):\xxQ  \to \xxQ,\qquad
 [L(t)\chi](s)  := S(t)\chi(s), \quad s\in (0,\ell),
\end{equation*}
and the mapping
\begin{equation*}
e:\xxQ \to \ltwob,
\qquad
\chi \mapsto \chi(\ell).
\end{equation*}
The solutions are understood in the sense of Theorem~\ref{teowell},
hence elements of $\xxQ$ posses additional regularity.
In particular, for any $\chi \in \xxQ$, one has
\begin{align} \label{reg-traj1}
\chi &\in L^{\infty}(0,\ell;\ltwob) \cap \Lbfi{2}(0,\ell;W^{1,2}(\Omega))
    \cap \Lbfi{p}(0,\ell;L^p(\Omega);
\\ \label{reg-traj2}
\chi_t & \in L^2(0,\ell;\wnegx) + L^{p'}(0,\ell;\lppx).
\end{align}
Also, thanks to Corollary~\ref{new:cor}, we can assume that
\begin{equation} \label{reg-traj3}
  \chi \in L^{\infty}(0,\ell;L^{q}_b(\Omega)).
\end{equation}
All the above estimates are independent of $\chi$ and $\barx$.
Consequently $\chi \in C([0,\ell];\ltwox)$ in the sense of representative,
and it thus makes sense to talk about point values of elements of $\xxQ$.
Continuity properties of the above introduced operators are 
summarized in the following:

\bete\label{te-cont1}
\begin{enumerate}

\item $S(t):\ltwox \to \ltwox $ are Lipschitz continuous
uniformly w.r.t.~$t\in [0,T]$;

\item $L(t):L^2(0,\ell;\ltwox) \to L^2(0,\ell;\ltwox)$ are 
Lipschitz continuous uniformly w.r.t.~$t\in [0,T]$;

\item $e:L^2(0,\ell;\ltwox) \to \ltwox$ is Lipschitz continuous.

\end{enumerate}

\ente
\begin{proof}
Let $u_1$, $u_2$ be weak solutions. Subtract the equations and
test by $w:=u_1-u_2$. We have
\begin{align}\no
&\frac12 \ddt \io |w(x,t)|^2 \dwx
\\
\no
&+
\io \big( a(\nabla u_1(x,t)) - a(\nabla u_2(x,t)) \big)
\cdot \big( \nabla w(x,t) - w(x,t)
\frac{x-\barx}{|x-\barx|} \big) \dwx
\\
&+
\io \big( f(u_1(x,t)) + \h(x,\nabla u_1(x,t)) - f(u_2(x,t)) -\h(x,\nabla u_2(x,t)) \big) w(x,t) \dwx =0.
\end{align}
Invoking \eqref{hpa1}-\eqref{hpa2}, \eqref{hpf3}
and \eqref{hph} and using Young's
inequality, one deduces
\[
\ddt \io |w(x,t)|^2 \dwx
+ \kappa \io |\nabla w(x,t)|^2 \dwx
\le c_1 \io  |w(x,t)|^2 \dwx.
\]
Integration over $t\in (t_1,t_2)$ yields
\begin{gather*}
 \io |w(x,t_2)|^2 \dwx
 + 
 \kappa \int_{\Omega \times (t_1,t_2)} 
 |\nabla w(x,t)|^2 \dwxt
 \\
 \le 
 \io |w(x,t_1)|^2 \dwx
 +
 c_1  \int_{\Omega \times (t_1,t_2)} 
  |w(x,t)|^2 \dwxt,
\end{gather*}
and Gronwall's lemma applied to 
\[
Y(t) := \io  |w(x,t)|^2 \dwx 
+ \kappa \int_{\Omega \times (t_1,t)} 
 |\nabla w(x,s)|^2 \dwx \dis
\]
implies the basic estimate
\begin{equation} \label{diffw1}
\sup_{t\in [t_1,t_2] }
\io |w(x,t)|^2 \dwx
+
\kappa \int_{\Omega \times (t_1,t_2)} |\nabla w(x,t)|^2 \dwxt
\le
c_2 \io |w(x,t_1)|^2 \dwx.
\end{equation}
Part 1 of the theorem follows immediately. One also has
\begin{equation*}
\norm{w(t+s)}{\ltwox}{2} \le c_3 \norm{w(s)}{\ltwox}{2},
\qquad
\norm{w(\ell)}{\ltwox}{2} \le c_4 \norm{w(s)}{\ltwox}{2},
\end{equation*}
for any $s\in(0,\ell)$, $t\in(0,T)$, where the constants 
$c_3$, $c_4$ only depend on $T$. Then, 
integrating the above relations over $s$ yields
parts 2 and 3 of the theorem, respectively.
\end{proof}
\beos
We can establish even stronger continuity of $S(t)$.
From the above theorem, one has
\[
\norm{w(t_2)}{\ltwox}{2} \le c_4 \norm{w(t_1)}{\ltwox}{2};
\]
multiplying by $\phi(\barx)$ and taking suprema over $\barx$,
together with Theorem~\ref{teoequiv1}, yields the continuity
of $S(t)$ with respect to the $L^2_{b,\phi}(\Omega)$-norm.
\eddos

The existence of a global attractor is now proved in a straightforward
manner. Recall that, following \cite{BaVi92}, a set $\atr$ is called
$(X,Y)$-attractor for the dynamical system $(S(t),X)$, provided that
$\atr$ is fully invariant, compact in the topology $Y$, and attracts
bounded subsets of $X$ uniformly in the topology of $Y$.

\bete\label{teoatrakt}
The dynamical system $(S(t),\ltwob)$ has a $(\ltwob,\ltwoloc)$-attractor.
\ente
\begin{proof}
1. We first establish the attractor for $(L(t),\xxQ)$.
Recalling Theorem~\ref{teoequiv4} above, it follows from
\eqref{reg-traj1}, \eqref{reg-traj2} that $\xxQ$ is bounded in
each of the seminorms
\begin{align*} 
\chi &\in L^2(0,\ell;W^{1,2}(C_k)) \cap L^p(0,\ell;L^p(C_k)),
\\ 
\chi_t & \in L^2(0,\ell;W^{-1,2}(C_k)) + L^{p'}(0,\ell;L^{p'}(C_k)).
\end{align*}
By the Aubin-Lions Lemma, we then have compactness in 
$L^2(0,\ell;L^2(C_k))$ for any $k$; invoking the boundedness of $\xxQ$ in
$L^{\infty}(0,\ell;\ltwob)$, we have indeed the compactness
in $L^2(0,\ell;\ltwox)$. Recalling the continuity of $L(t)$, we
deduce the existence of $\atrQ$, the global attractor for $(L(t),\xxQ)$,
by standard arguments.

2. Set
\begin{equation} \label{a-atr}
\atr := e(\atrQ).
\end{equation}
From the continuity of $e$ and the equivalence \eqref{loc-x}
one immediately obtains that
this is the desired $(\ltwob,\ltwoloc)$-attractor for $(S(t),\ltwob)$.
\end{proof}

\beos\label{exglobattr}
The  \emph{existence} of the global
attractor can be proven solely in virtue of the regularity
established in Theorem~\ref{teowell} above. Of course, in this case
we can no longer choose  $\bnul$ bounded in $L^{q}_b(\Omega)$.
Moreover, extensions to systems are possible on account of Remarks~\ref{vector1}
and \ref{vector2}.
\eddos

\section{Entropy estimates}
\label{entroest}

The aim of the last section is to study finite-dimensionality
of the attractor. 
As is well known, for dissipative equations in the case of a
\emph{bounded} domain $\Omega$, the attractor $\atr_\Omega$ satisfies
\begin{equation} \label{est-f}
\hent{\varepsilon}{\atr_\Omega}{L^2(\Omega)}
\le c_0 \vol(\Omega) \ln\frac1\varepsilon,
\qquad \eps \in (0,\epsn).
\end{equation}
Here the constant $c_0$
only depends on the structural properties of the equation, but not
on the size of $\Omega$. In particular, we have finite fractal dimension
of $\atr_\Omega$. Such an estimate being meaningless if $\Omega$ has infinite volume,
we will follow \cite{Ze01,Ze03} to estimate the entropy
of $\atr$ in the seminorm $L^2_b(\oo)$, where $\oo$ is a suitable
bounded subdomain of $\Omega$.

Our main result is the following theorem.
\bete\label{teo-main}
Let $d\leq 3$ and set
\[
\Omxr := \Omega \cap B_{R}(x_0,\RR^d) = B_{R}(x_0,\RR^d).
\]
Then, there exist $c_0$, $c_1$ and $\epsn>0$, such that, 
for any $x_0 \in \Omega$, $R\ge 1$ and $\varepsilon \in (0,\epsn)$ one has
\begin{equation} \label{est-main}
\hent{\varepsilon}{\atr}{L^2_{b}(\Omxr)}
\le c_0 \left( R + c_1 \ln\frac1\varepsilon \right)^d \ln\frac1\varepsilon.
\end{equation}
\ente
The rest of this section is devoted to the proof of this result. Remark 
that \eqref{est-main} is completely analogous to \eqref{est-f},
but for the ``extra term'' $c_1 \ln\frac1\varepsilon$.
Heuristically, the finer description of $\atr$ one seeks, the
larger portion of $\Omega$ influences the dynamics. Moreover,
the optimality of this estimate is suggested by the results
of \cite{Ze03} where a similar bound is proved to be sharp,
albeit in a different regularity setting.
\par
Given $x_0 \in \Omega$ and $R\ge 1$, we set
\begin{equation} \label{def-psi}
\Psixr:=
\begin{cases}
1; & |x-x_0| \le R + \sqrt{d},
\\
\exp\big((R+\sqrt{d} - |x-x_0|)/2\big); & \textrm{otherwise.}
\end{cases}
\end{equation}
Clearly, one has
\[
\hent{\varepsilon}{\atr}{L^2_{b}(\Omxr)} 
\le
\hent{\varepsilon}{\atr}{L^2_{b,\Psixr}(\Omega)},
\]
hence it is enough to estimate the right-hand side. As usual, one
arrives at such a result through suitable iterative coverings obtained
by combining the ``smoothing property'' of solution operators with
compact embeddings in the appropriate function spaces. As in the
previous section, we will rely on the natural parabolic estimates.
Let us start by an improved continuity result for the evolution
operators.

\bete\label{te-cont2}
Let $\phi$ be an admissible weight function of growth rate $\mu <1$.
Then,
\begin{enumerate}

\item $L(t):\XQbfi \to \XQbfi$ are Lipschitz continuous
uniformly w.r.t. $t\in [0,T]$;

\item $e:\XQbfi \to \ltwobfi$ is Lipschitz continuous.

\end{enumerate}

\ente
\begin{proof}
It follows from \eqref{diffw1} that
\[
\io |w(x,t+s)|^2 \dwx \le c_4 \io |w(x,s)|^2 \dwx,
\qquad t\in [0,T],
\]
where $c_4$ only depends on $T$. Hence,
integrating in $\dis$ over $(0,\ell)$,
\[
\int_{\Omega \times (t,t+\ell)} |w(x,s)|^2 \dwx \dis
\le
c_4
\int_{\Omega \times (0,\ell)} |w(x,s)|^2 \dwx \dis.
\]
Applying $\sup_{\barx \in \Omega} \phi(\barx)$ and using
Theorem~\ref{teoequiv4} again, one obtains -- in terms of trajectories --
\begin{equation} \label{diff-Lt}
\norm{L(t)\chi_1 - L(t)\chi_2}{\XQbfi}{}
\le c_5
\norm{\chi_1 - \chi_2}{\XQbfi}{}.
\end{equation}
This proves part 1. Analogously, one deduces
\[
\io |w(x,\ell)|^2 \dwx
\le c_6
\iq |w(x,t)|^2 \dwxt,
\]
and thus
\begin{equation} \label{diff-e}
\norm{e(\chi_1) - e(\chi_2)}{\ltwobfi}{}
\le c_5
\norm{\chi_1 - \chi_2}{\XQbfi}{},
\end{equation}
i.e., part 2.
\end{proof}

A key step towards our entropy estimate is the following 
``smoothing property''
of the operator $L(t)$ -- a sort of typical result in the spirit of
the method of trajectories.

\bete\label{te-smoo}
Let $\phi$ be an admissible weight function of growth rate $\mu <1$.
Then for any $\chi_1$, $\chi_2 \in \xxQ$, one has
\begin{align}   \label{sm1}
\norm{L(\ell)\chi_1 - L(\ell)\chi_2}{\VQbfi}{}
    &\le K_1 \norm{\chi_1 - \chi_2}{\XQbfi}{},
\\  \label{sm2}
\norm{\dt(L(\ell)\chi_1 - L(\ell)\chi_2)}{\ZQbfi}{}
    &\le K_2 \norm{\chi_1 - \chi_2}{\XQbfi}{}.
\end{align}
where $K_1$, $K_2$ only depend on 
the constants $\mu$ and $c$ characterizing
the growth of $\phi$ in \eqref{weight-gr}.
\ente
\begin{proof} Let $u$, $v$ be the weak solutions such that
$u|_{[0,\ell]}=\chi_1$, $v|_{[0,\ell]}=\chi_2$, and set $w:=u-v$.

\noindent
STEP 1. One deduces from \eqref{diffw1} that
\[
\int_{\Omega \times (\ell,2\ell)}
\big( |\nabla w(x,t)|^2 + |w(x,t)|^2 \big) \dwxt
\le
c_1 \io |w(x,s)|^2 \dwx,
\qquad \forall s\in(0,\ell).
\]
Integrating over $s$ and applying $\sup_{\barx}\phi(\barx)$,
in view of Theorems~\ref{teoequiv1}, \ref{teoequiv4},
yields \eqref{sm1}.

\noindent STEP 2.
We will show that
\[
\norm{\dt w}{\ZQbfi}{} \le c_2 \norm{w}{\VQbfi}{},
\]
which combined with \eqref{sm1} implies \eqref{sm2}.
Using equation \eqref{eq} and Theorem~\ref{teoequiv4}, we can write
\begin{gather}\no
\norm{\dt w}{\ZQbfi}{} 
\le c_3  \sup_z \sup_{\barx}   \phi(\barx)
\int_{\Omega \times (0,\ell)} \dt w z \dwx
\\
 \no
= \sup_z \sup_{\barx} \phi(\barx)\int_{\Omega \times (0,\ell)}
\Big[\big( a(\nabla u) - a(\nabla v) \big)\cdot
\big(\nabla z - z \frac{x-\barx}{|x-\barx|} \big)
+ \big( \h(x,\nabla u) - \h(x,\nabla v) \big) z \\
\label{new:11}
+ \big( f(u) - f(v) \big)z\Big] \dwxt.
\end{gather}
The first supremum is taken over $z\in \VQbfi$
with unit norm.
Recalling \eqref{hpa2}, the first two terms on the right-hand side get estimated as
\begin{gather*}
c_4 \phi(\barx) \iq
|\nabla w| \big( |\nabla z| + |z| \big) \dwxt
\\
\le
c_4 \Big( \phi(\barx) \iq |\nabla w|^2 \dwxt \Big)^{1/2}
 \Big( \phi(\barx) \iq (|\nabla z|+ |z|)^2 \dwxt \Big)^{1/2}
 \\
\le c_5
\norm{w}{\VQbfi}{} \norm{z}{\VQbfi}{}
= c_5 \norm{w}{\VQbfi}{}.
\end{gather*} 
Invoking \eqref{hpf2}, the last term 
in the \rhs\ of \eqref{new:11}
is estimated as
\begin{align} \no
  & c_6 \phi(\barx) \iq \big(1+|u|+|v|)^{p-2} |w| |z| \dwxt\\
 \label{kn1}  
  & \mbox{}~~~~~~~~~~  
   \le c_7 \sup_k \phi(x_k) \int_0^\ell 
     \norm{ (1+|u|+|v|)^{p-2} wz }{L^1(C_k)}{} \,\dit.
\end{align}
We have used \eqref{ekvxp} with $p=1$. 
Furthermore, thanks to the embedding $H^1(\Omega)\hookrightarrow L^6(\Omega)$
and the additional regularity 
\eqref{reg-traj3} 
(where we take $q=3(p-2)/2$),
we can estimate 
by H\"older's inequality
\[
  \norm{(1+|u|+|v|)^{p-2}wz}{1}{}
   \le \norm{(1+|u|+|v|)}{\frac{3(p-2)}{2}}{p-2}
   \norm{w}{6}{} \norm{z}{6}{}
  \le c_8 \norm{w}{W^{1,2}(C_k)}{}
  \norm{z}{W^{1,2}(C_k)}{}.
\]
Hence, \eqref{kn1} can further be estimated as
\[
  c_9 \Big( \sup_{k} \phi(x_k) \int_0^\ell \norm{w}{W^{1,2}(C_k)}{2} \, \dit \Big)^{1/2}
  \Big( \sup_{k} \phi(x_k) \int_0^\ell \norm{z}{W^{1,2}(C_k)}{2} \, \dit \Big)^{1/2}
  \le c_{10} \norm{w}{\VQbfi}{}.
\]
This finishes the proof.
\end{proof}
What we have just shown is the Lipschitz continuity of $L(\ell)$ from
$\XQbfi$ into $\WQbfi$, where the latter space was defined in 
\eqref{W-norm}.
However -- and  this is the peculiar feature
of the analysis in unbounded domains -- the space $\WQbfi$ is NOT
compactly embedded into $\XQbfi$. The compactness can only be employed
using seminorms related to restrictions to bounded sets $\oo \subset \Omega$
(cf.\ Lemma~\ref{le-aub_res}) which also exhibit the correct dependence
on the volume of the domain of restriction $\oo$.

Last ingredient is to employ the boundedness of $\xxQ$ in $\ltwob$
together with the decay of the weight $\Psixr$ to localize the entropy
of an attractor to a bounded domain up to some error. This estimate
is actually the origin of the ``extra term'' in \eqref{est-main}.

\bele \label{lm-decay}
Let $\epsn>0$ be given. Then there exists $c_1$ such that,
for any $x_0 \in \RR^d$, $R\ge1$ and $\eps \in (0,\epsn)$,
having set
\[
\Reps:= R + c_1\left( 1 + \ln\frac1\varepsilon \right),
\]
for arbitrary $\chi_1$, $\chi_2 \in \xxQ$, one has
\[
\norm{\chi_1 - \chi_2}{\XQbpsi}{}
\le 
\max\big\{ 
\norm{\chi_1-\chi_2}{L^2_{b,\Psixr}(0,\ell; L^2(\Omega_{x_0,\Reps}))}{},
\varepsilon
    \big\}.
\]
\enle
\begin{proof}
Recall that
\begin{align*}
  \norm{\chi_1 - \chi_2}{\XQbpsi}{2}
  = \max\big\{ 
  &\sup_{x_k \notin \Omega_{x_0,\Reps}}
   \Psixr(x_k) \norm{\chi_1 -\chi_2}{L^2(0,\ell;L^2(C_k))}{2},
\\
&\norm{\chi_1-\chi_2}{L^2_{b,\Psixr}(0,\ell; L^2(\Omega_{x_0,\Reps}))}{}
\big\}.
\end{align*}
However, thanks to the decay of $\Psixr$ and the boundedness of $\xxQ$,
the first term is automatically smaller than $\varepsilon^2$ due to
proper choice of constant $c_1$.
\end{proof}

\medskip \noindent
We are now ready to prove the main result.

\medskip \noindent
{\sc Proof of Theorem~\ref{teo-main} }
In what follows, $\Psixr$ is the weight function defined in \eqref{def-psi}.
Remark that it has growth rate $\mu=1/2$, and satisfies \eqref{weight-gr}
with $c=1$ independently of and $x_0\in \RR^d$, $R\ge1$.

\noindent
1. First (and the key) step of the proof is the recurrent estimate
\begin{equation} \label{recc}
\hent{\alpha/2}{\atrQ}{\XQbpsi} 
\le
\hent{\alpha}{\atrQ}{\XQbpsi} + c_0 \left( R + c_1 \ln\frac1\alpha \right)^d
.
\end{equation}
Indeed, let
\[
\atrQ \subset \bigcup_m B_{\alpha}(\chi_m;\XQbpsi)
.
\]
Thanks to Theorem~\ref{te-smoo} and invariance of $\atrQ$, we then deduce
that, for some $\tilde{\chi}_m \in \calX$ and some $\kappa>0$
\[
\atrQ \subset \bigcup_m B_{\kappa \alpha}(\tilde{\chi}_m; W_{b,\psi}(Q))
.
\]
By Lemma~\ref{le-aub_res}, each of the latter balls can be covered
so that
\[
\hent{\alpha/2}{B_{\kappa \alpha}(\tilde{\chi}_m; W_{b,\psi}(Q))}%
{X_{b,\psi}(\Omega_{x_0,R(\alpha/2)})}
\le
c_0  \left( R + c_1 \ln\frac1\alpha \right)^d
.
\]
Finally, by Lemma~\ref{lm-decay}, covering of $\atrQ$ 
by $\alpha/2$-balls in the
$X_{b,\psi}(\Omega_{x_0,R(\alpha/2)})$ seminorm is also
covering by $\alpha/2$-balls in the norm $\XQbpsi$.

\noindent
2. Choose $\epsn>0$ such that $\hent{\epsn}{\atrQ}{\XQbpsi} = 0$.
Given $\eps\in(0,\epsn)$, one picks $k\in\NN$ such that
\[
2^{-k} \epsn \le \eps < 2^{-k+1}\epsn
.
\]
Note that this means $k\le c\ln1/\eps$, 
at least provided $\eps$ is small enough. 
Then, using \eqref{recc}, one can estimate
\begin{align*}
  & \hent{\eps}{\atrQ}{\XQbpsi} \le \hent{2^{-k}\epsn}{\atrQ}{\XQbpsi}\\
  & \mbox{}~~~~~~~~~~ \le \sum_{l=1}^{k}
   \hent{2^{-l}\epsn}{\atrQ}{\XQbpsi} - \hent{2^{-l+1}\epsn}{\atrQ}{\XQbpsi} \\
  & \mbox{}~~~~~~~~~~ \le \sum_{l=1}^{k}
   c_0 \left( R + c_1 \ln\frac{2^{l-1}}{\epsn} \right)^d \\
  & \mbox{}~~~~~~~~~~
   \le c_0 \left( R + c_1 \ln\frac{1}{\eps} \right)^d \ln\frac1\eps.
\end{align*}

3. Finally, in view of Lipschitz continuity of $e$ (Theorem~\ref{te-cont2})
and \eqref{a-atr} we conclude
\[
\hent{\eps}{\atr}{L^2_{b,\Psixr}(\Omega)}
\le
\hent{\eps/\kappa}{\atrQ}{\XQbpsi}
,
\]
%
where $\kappa$ is the Lipschitz constant of the mapping $e$,
which is independent of the particular weight function $\Psixr$.
This finishes the proof.
\qed

\beos\label{vector3}
Recalling Remark~\ref{vector2}, we point out that Theorem~\ref{teo-main} can be
extended to systems where the diffusion operator is the vector Laplacian
$-\mathbf{\Delta}$ and the nonlinear convective term $h$ is replaced by 
$\mathbf{h}:\Omega\times\RR^{m\times d}\to\RR^m$ 
satisfying suitable generalizations of \eqref{hph} and \eqref{hph2}. 
\eddos




\vspace{5mm}

\noindent%
{\bf First author's address:}\\[1mm]
Maurizio Grasselli\\
Dipartimento di Matematica, Politecnico di Milano\\
Via E. Bonardi, 9,~~I-20133 Milano,~~Italy\\
E-mail:~~{\tt maurizio.grasselli@polimi.it}

\vspace{4mm}

\noindent%
{\bf Second author's address:}\\[1mm]
Dalibor Pra\v z\'ak\\
Katedra matematick\'e anal\'yzy\\
Matematicko-fyzik\'aln\'{\i} fakulta Univerzity Karlovy\\
Sokolovsk\'a 83, 186 75, Czech Republic\\
E-mail:~~{\tt prazak@karlin.mff.cuni.cz}

\vspace{4mm}

\noindent%
{\bf Third author's address:}\\[1mm]
Giulio Schimperna\\
Dipartimento di Matematica, Universit\`a degli Studi di Pavia\\
Via Ferrata, 1,~~I-27100 Pavia,~~Italy\\
E-mail:~~{\tt giusch04@unipv.it}


\begin{thebibliography}{99}

\bibitem{A}
F.~Abergel, {\sl Existence and finite dimensionality of the global attractor for evolution equations on
unbounded domains}, J.\ Differential Equations, {\bf 83} (1990), 85–-108.

\bibitem{ACDR1}
J.M.~Arrieta, J.W.~Cholewa, T.~Dlotko, and A.~Rodr\'{\i}guez-Bernal, 
{\sl Asymptotic behavior and
attractors for reaction diffusion equations in unbounded domains}, 
Nonlinear Anal., {\bf 56} (2004), 515-–554.

\bibitem{ACDR2}
J.M.~Arrieta, J.W.~Cholewa, T.~Dlotko, and A.~Rodr\'{\i}guez-Bernal,  
{\sl Dissipative parabolic equations in locally uniform spaces},
Math.\ Nachr., {\bf 280} (2007), 1643--1663. 

\bibitem{BaVi90}
A.V.~Babin and M.I.~Vishik, {\sl Attractors of partial differential evolution equations in unbounded
domain}, Proc.\ Roy.\ Soc.\ Edinburgh Sect.\ A, {\bf 116} (1990), 
221--243.

\bibitem{BaVi92}
A.V.~Babin and M.I.~Vishik,
``Attractors of Evolution Equations'',
North-Holland, Amsterdam, 1992.

\bibitem{Ba}
 V.~Barbu,
 ``Nonlinear Semigroups and Differential Equations in Banach Spaces'',
 Noord\-hoff,
 Leyden,
 1976.
  
\bibitem{Br}
 H.~Br\'ezis,
 ``Op\'erateurs Maximaux Monotones et Semi-groupes de Contractions
   dans les Espaces de Hilbert'',
 North-Holland Math.\ Studies {\bf 5},
 North-Holland,
 Amsterdam,
 1973.
 
\bibitem{EM}
M.~Efendiev and A.~Miranville, {\sl 
Finite-dimensional attractors for reaction-diffusion equations in $\RR^n$ 
with a strong nonlinearity},
Discrete Contin.\ Dynam.\ Systems, {\bf 5} (1999), 399--424.  
  
\bibitem{EZ}
 M.A.~Efendiev and S.V.~Zelik, 
 {\sl The attractor for a nonlinear reaction-diffusion system in 
  an unbounded domain},
 Comm.\ Pure Appl.\ Math.,
 {\bf 54} (2001), 
 625--688.  
 
\bibitem{EZ2}
M.A.~Efendiev and S.V.~Zelik, 
{\sl Upper and lower bounds for the Kolmogorov entropy of the
attractor for the RDE in an unbounded domain}, 
J.\ Dynam.\ Differential Equations, {\bf 14} (2002), 369–-403.
 
\bibitem{FLST}
E.~Feireisl, Ph.~Lauren\c{c}ot, F.~Simondon, and H.~Tour\'e, 
{\sl Compact attractors for reaction-diffusion
equations in $\RR^n$}, 
C.\ R.\ Acad.\ Sci.\ Paris S\'er.\ I Math., {\bf 319} (1994), 147–-151.
 
\bibitem{Fe}
 E.~Feireisl, ``Dynamics of Viscous Compressible Fluids'',
 Oxford Lecture Series in Mathematics and its Applications, 26. 
 Oxford University Press, Oxford, 2004. 
 
\bibitem{GP}
M.~Grasselli and D.~Pra\v z\'ak, {\sl Exponential attractors for a class of reaction-diffusion 
problems with time delays}, J.\ Evol.\ Equ., {\bf 7} (2007), 649--667.
 
\bibitem{Li}
 J.-L.~Lions,
 ``Quelques M\'ethodes de R\'esolution des Probl\`emes aux Limites
    non Lin\'eaires'' (French),
 Dunod, Gauthier-Villars,
 Paris, 1969.
  
\bibitem{MP}
 J.~M\'alek and D.~Pra\v z\'ak, 
 {\sl Large time behavior via the method of $l$-trajectories},
 J.~Differential Equations,
 {\bf 181} (2002), 
 243--279. 
 
\bibitem{Me}
S.~Merino, {\sl On the existence of the compact global attractor for semilinear reaction diffusion
systems on $\RR^N$}, J.\ Differential Equations, {\bf 132} (1996), 87-–106.

\bibitem{MZ}
A. ~Miranville and S.~ Zelik, {\sl Attractors for
    dissipative partial
    differential equations in bounded and unbounded domains},
    Evolutionary equations. Vol. IV, 103--200, Handb. Differ. Equ.,
    Elsevier/North-Holland, Amsterdam 2008.

\bibitem{Pr}
M.~Prizzi, {\sl A remark on reaction-diffusion equations in unbounded domains}, 
Discrete Contin.\ Dyn.\ Syst., {\bf 9} (2003), 281--286. 

\bibitem{PR}
M.~Prizzi and K.P.~Rybakowski, {\sl Attractors for reaction-diffusion equations on arbitrary unbounded domains}, 
Topol.\ Methods Nonlinear Anal., {\bf 30} (2007), 251--277. 
 
   
    
\bibitem{Si}
 J.~Simon,
 {\sl Compact sets in the space {$L^p(0,T;B)$}},
 Ann.\ Mat.\ Pura Appl.\ (4),
 {\bf 146} (1987),
 65--96.

\bibitem{Te}
 R.~Temam,
 ``Infinite-Dimensional Dynamical Systems in Mechanics and Physics'',
 Second Edition,
 Springer-Verlag, New York, 1997.
 
\bibitem{Wa}
B.~Wang, 
{\sl Attractors for reaction-diffusion equations in unbounded domains}, 
Phys.\ D, {\bf 179} (1999), 41–-52. 

\bibitem{Ze01}
 S.V.~Zelik,
 {\sl The attractor for a nonlinear reaction-diffusion system in the
  unbounded domain and Kolmogorov's $\epsilon$-entropy},
 Math.\ Nachr.,
 {\bf 232} (2001), 
 129--179.

\bibitem{Ze03}
 S.V.~Zelik, 
 {\sl Attractors of reactions-diffusion systems in unbounded domains
  and their spatial complexity},
 Comm.\ Pure Appl.\ Math.,
 {\bf 56} (2003), no. 5, 
 584--637.  
 


\end{thebibliography}
\end{document}